\documentclass[11pt,a4paper]{article}
\usepackage[T1]{fontenc}
\usepackage{lmodern}
\usepackage[margin=27mm]{geometry}
\usepackage{amsmath,amssymb,amsthm,mathtools}
\usepackage{microtype}
\usepackage{enumitem}
\usepackage{tikz-cd}
\usepackage{cite}
\usepackage{xcolor}
\usepackage{hyperref}
\hypersetup{colorlinks=true,linkcolor=blue!45!black,citecolor=blue!45!black,urlcolor=blue!45!black,pdftitle={The Maximum Length of MDS and NMDS Algebraic Geometry Codes}}
\numberwithin{equation}{section}

\theoremstyle{plain}
\newtheorem{theorem}{Theorem}[section]
\newtheorem{lemma}[theorem]{Lemma}
\newtheorem{proposition}[theorem]{Proposition}
\newtheorem{corollary}[theorem]{Corollary}

\theoremstyle{definition}
\newtheorem{definition}[theorem]{Definition}

\theoremstyle{remark}
\newtheorem{remark}[theorem]{Remark}

\newcommand{\F}{\mathbb F}
\newcommand{\wh}{\widehat}

\newcommand{\supp}{\operatorname{Supp}}

\newcommand{\CL}{C_{\mathcal L}}

\title{The Maximum Length of MDS and NMDS Algebraic Geometry Codes}
\author{}
\date{}

\begin{document}

\begin{center}
    {\Large On the Maximum Length of MDS and NMDS Algebraic Geometry Codes}

    \vspace{1em}

    {\large Lian Duan\footnotemark[1] 
    and Dongchun Han\footnotemark[2]}
\end{center}

\footnotetext[1]{Institute of Mathematical Sciences,
ShanghaiTech University, Shanghai 201210, P.R. China.
Email: \texttt{duanlian@shanghaitech.edu.cn}}

\footnotetext[2]{School of Mathematics,
Southwest Jiaotong University, Chengdu 610000, P.R. China.
Email: \texttt{handongchun@swjtu.edu.cn}}

\vspace{0.6em}

\begin{abstract}
We study the maximum length of MDS and near-MDS (NMDS) algebraic geometry codes over $\mathbb{F}_q$. For sufficiently large dimension and codimension, we prove that MDS codes of genus $g\ge 2$ and NMDS codes of genus $g\ge 3$ have length less than $20g\sqrt{q}$, and obtain the sharp upper bound $\#X(\mathbb{F}_q)/2+\sqrt{q}$ for NMDS codes on a genus-two curve $X$. For each fixed genus $g\ge2$ and infinitely many $q$, Hu and Ma~\cite{HuMaConstruction} constructed MDS codes of length $\sqrt{q}+2gq^{1/4}+1<2\sqrt{q}$, and we construct NMDS codes of the same length. Thus the square-root upper bounds have the correct order of magnitude. We also show that, for each fixed genus, the MDS and NMDS conjectures hold for algebraic geometry codes over all sufficiently large finite fields.
\end{abstract}
	
	\section{Introduction}\label{Sect: intro}

	Let $C$ be a linear $[n,k,d]_q$ code. The Singleton bound gives $d\le n-k+1$. The code is \emph{maximum distance separable} (MDS) if equality holds, and \emph{near-MDS} (NMDS) if both $C$ and its dual have Singleton defect one, that is, $d(C)=n-k$ and $d(C^\perp)=k$. We call $C$ nontrivial if $2\le k\le n-2$. MDS and NMDS codes are central objects in coding theory and finite geometry. They attain or nearly attain the largest minimum distance allowed by the Singleton bound. Geometrically, MDS codes correspond to arcs in projective spaces, and three-dimensional NMDS codes correspond to $(n,3)$-arcs in the projective plane \cite{MS,BL,DL,LR}.

	The \emph{MDS conjecture} states that a nontrivial linear MDS code has length at most $q+1$, with $q+2$ allowed when $q$ is even and $k\in\{3,q-1\}$. The geometric approach goes back to Segre's work on ovals and was developed by Bruen, Thas, and Blokhuis and by Storme and Thas \cite{Segre,BTB,ST}. Ball settled the conjecture over prime fields, and Ball and De Beule extended the known range over extension fields \cite{Ball,BDB}. We refer to Ball and Lavrauw \cite{BL} for further developments. More recently, Fan obtained a range of dimensions linear in $q$ in odd characteristic, while Kumar and, independently, Bocong Chen and Jing Huang gave proofs in even characteristic \cite{Fan,Kumar,CHMDS}.

	The \emph{NMDS conjecture} predicts that a nontrivial NMDS code has length at most $2q+2$ for $q\ge5$ \cite[Conjecture~5.2]{LR}. In dimension three, Thas's work on planar arcs gives the stronger bound $n\le2q+1$ for $q>3$ \cite{Thas1975,LR}. Beyond dimension three, exact results are known for some small fields and particular dimensions, and the conjecture remains open in general \cite{DL,DLsmall,LR}.

	Goppa's construction gives a natural class of codes in which to study these problems. Let $X$ be a smooth, projective, geometrically integral curve of genus $g$ over $\F_q$. Let $D=P_1+\cdots+P_n$ be a sum of distinct rational points of $X$, and let $G$ be a rational divisor whose support is disjoint from that of $D$. The algebraic geometry code $\CL(D,G)$ is obtained by evaluating the Riemann--Roch space $L(G)$ at $P_1,\ldots,P_n$. In genus zero this gives generalized Reed--Solomon codes \cite{Goppa,TV,Liu,Stichtenoth}. Tsfasman, Vl\u{a}du\c{t}, and Zink showed that curves with many rational points yield families of codes surpassing the Gilbert--Varshamov bound \cite{TVZ}. Garc\'ia and Stichtenoth later gave an explicit tower attaining the Drinfeld--Vl\u{a}du\c{t} bound \cite{GS}.

	For elliptic codes, the group law relates the minimum distance to sums of distinct rational points. This connection appears in early work of Katsman and Tsfasman, Shokrollahi, and Munuera \cite{KT,Shokrollahi,MunueraEll,Munuera}. Li, Wan, and Zhang conjectured that the length of a nontrivial MDS elliptic code is asymptotically at most half the number of rational points on the curve \cite{LWZ}. Han and Ren first proved the conjecture in a range of dimensions and then settled it in full \cite{HRTight,HR}. More recently, Bocong Chen and Jing Huang sharpened the half-point bound, and Han and Ren obtained stronger bounds for curves with an odd number of rational points \cite{CHRestricted,HRCover}.

Earlier work on MDS algebraic geometry codes of higher genus studied the geometric MDS conjecture~\cite{Munuera,HaoChen,Walker,CY,deBoer}. Zhu, Chen, Peng, Zhang, and Zhao~\cite{ZhuEtAl} constructed long MDS codes on maximal curves, and Hu and Ma~\cite{HuMaConstruction} constructed MDS codes of square-root length for every fixed genus $g\ge2$. Haojie Chen~\cite{Chen2} proved an upper bound of $2q/3+O(\sqrt q)$ for genus-two MDS codes. For any fixed positive genus, Chen, Hu, and Ma~\cite{CHM} obtained an upper bound of $(1/2+o(1))q$ when both the dimension and codimension grow faster than $\log q$. For NMDS algebraic geometry codes, Zhu et al.~\cite{ZhuEtAl} gave
explicit constructions on genus-two curves with split Jacobians. For fixed genus, how large can MDS and NMDS algebraic geometry codes
be when both the dimension and codimension are large?

	We prove square-root upper bounds for MDS codes of genus at least two and NMDS codes of genus at least three when the dimension and codimension are sufficiently large. For every fixed genus $g\ge2$, we construct NMDS codes of length $\sqrt q+2gq^{1/4}+1$ over infinitely many square fields. Together with the MDS construction of Hu and Ma~\cite{HuMaConstruction} at the same length, this shows that the square-root bounds have the correct order of magnitude. In genus two, we give a sharp upper bound for NMDS codes in terms of the number of rational points on the curve.

	\subsection*{Main results}

	The dual of an algebraic geometry code on $X$ is again an algebraic geometry code on $X$, and duality preserves the MDS and NMDS properties (see Lemmas~\ref{Lem: dual_AG},~\ref{Lem: dual_MDS} and Remark~\ref{Lem: dual_NMDS} below). Also, we focus on the nontrivial codes. Therefore we may assume $2\le k\le n/2$. All logarithms are natural.
For $g\ge2$, set \(L_g(q)=2(g+3)\log q+4\log\!\left(16\binom{2g-2}{g-1}\right)\). For $g\ge3$, set \(L'_g(q)=1+2(g+4)\log q+4\log\!\left(16\binom{2g-2}{g-2}\right)\).

\begin{theorem}\label{thm:length}
Let $C$ be an MDS algebraic geometry code of genus $g\ge 2$,
length $n$, and dimension $g+2\le k\le n/2$. Then
\begin{equation}\label{eq:mds-length-bounds}
n\le
\begin{cases}
65g\,q^{\frac{g+1}{k}},  & g+2\le k<2g,\\
33g\,q^{\frac{g+1}{2g}}, & 2g\le k<L_g(q),\\
20(g-1)\sqrt{q},         & k\ge L_g(q).
\end{cases}
\end{equation}
\end{theorem}

For fixed genus, the preceding bounds yield the following asymptotic consequences as $q\to\infty$.
	\begin{corollary}\label{thm:mds}
		Fix $g\ge2$ and $\varepsilon>0$. For all $q$ sufficiently large in terms of $g$ and $\varepsilon$, every nontrivial MDS algebraic geometry code of genus $g$ with $k\le n/2$ satisfies
		\[
			n\le
			\begin{cases}
				q+1,& k<L_g(q),\\
				q^{1/2+\varepsilon},& k\ge L_g(q),
			\end{cases}
		\]
		with $q+2$ allowed in the first case when $q$ is even and $k=3$. In particular, for every fixed genus, the MDS conjecture holds for algebraic geometry codes over all sufficiently large finite fields.
	\end{corollary}

\begin{theorem}\label{thm:nmds-length}
Let $C$ be an NMDS algebraic geometry code of genus $g\ge 3$,
length $n$, and dimension $g+2\le k\le n/2$. Then
\begin{equation}\label{eq:nmds-length-bounds}
n\le
\begin{cases}
65g\,q^{\frac{g+2}{k+1}}, & g+2\le k<2g-1,\\
33g\,q^{\frac{g+2}{2g}},  & 2g-1\le k<L'_g(q),\\
20(g-1)\sqrt{q},          & k\ge L'_g(q).
\end{cases}
\end{equation}
\end{theorem}

For fixed genus $g$ and sufficiently large $q$, the NMDS conjecture for algebraic geometry codes follows directly from the Hasse--Weil bound. When $g\ge 3$ and $k\ge L'_g(q)$, our theorem gives a sharper upper bound of order $\sqrt{q}$.
	\begin{corollary}\label{cor:nmds}
		Fix $g\ge3$ and $\varepsilon>0$. For all $q$ sufficiently large in terms of $g$ and $\varepsilon$, every nontrivial NMDS algebraic geometry code of genus $g$  with $k\le n/2$ satisfies
		\[
			n\le
			\begin{cases}
				2q+2,& k<L'_g(q),\\
				q^{1/2+\varepsilon},& k\ge L'_g(q).
			\end{cases}
		\]
	\end{corollary}

\begin{remark}
    Hu and Ma~\cite[Theorem~1 and Corollary~3]{HuMaConstruction} constructed, for every fixed genus $g\ge2$ and infinitely many $q$, MDS algebraic geometry codes of length $n=\sqrt q+2gq^{1/4}+1$ for every dimension $g+1\le k\le n-g$. In Section~\ref{Sect: sqrt_root_length}, we construct NMDS codes with the same length and dimension range. Thus the exponent $1/2$ is optimal in the large-dimension upper bounds for MDS codes of genus $g\ge2$ and NMDS codes of genus $g\ge3$.
\end{remark}

The behavior of NMDS algebraic geometry codes in genus two differs from that in higher genus. For every fixed genus $g\ge3$, Theorem~\ref{thm:nmds-length} gives a square-root upper bound when the dimension is sufficiently large, whereas the following proposition shows that NMDS codes of genus two can have length proportional to $q$. We restrict our discussion of this phenomenon to genus-two NMDS codes. We expect a similar extension for codes whose Singleton defects, as well as those of their duals, equal an integer $\delta\ge1$, with the critical genus being $g=\delta+1$.

\begin{proposition}\label{prop:nmds-genus-two}
Let $C$ be an NMDS algebraic geometry code on a genus-two curve $X$ over $\F_q$, of length $n$ and dimension $k$. If $q\ge2500$ and $60\log q\le k\le n/2$, then
\begin{equation}\label{eq:n-genus2}
n\le\frac{\#X(\F_q)}{2}+\sqrt q.
\end{equation}
Moreover, for infinitely many $q$, there exist genus-two curves and NMDS codes attaining equality in \eqref{eq:n-genus2} for every dimension $2\le k\le n-2$.
\end{proposition}

Our approach builds on the work of Han and Ren~\cite{HRTight}, who used the group law on an elliptic curve to bound the length of MDS elliptic codes. Their argument uses the fact that an elliptic curve is isomorphic to its Jacobian. For a curve of fixed genus $g\ge 2$ over a large finite field, the rational points of the curve form only a small fraction of the rational points of its Jacobian. Chen, Hu, and Ma addressed this difficulty using character sums on the Jacobian. They expressed the MDS condition in terms of the vanishing of certain Riemann--Roch spaces~\cite[Proposition~5]{CHM}, which led them to study sums of distinct evaluation points in the Jacobian. By comparing the contribution of the trivial character with the total contribution of the nontrivial characters, they obtained an upper bound on the code length~\cite[Proposition~10]{CHM}. Some of the ideas in the present paper also grew out of joint work by Xia, Ren, and Han~\cite{XRH} on the MDS and NMDS conjectures for algebraic geometry codes.

Our proof follows the same general strategy, but develops a different approach to estimating the contributions of nontrivial characters. Chen, Hu, and Ma~\cite{CHM} control the contributions of nontrivial characters by using a uniform upper bound for the absolute values of the evaluation sum for every nontrivial character. Instead, we divide the characters into two classes according to whether the absolute values of these sums exceed a prescribed threshold. For characters below the threshold, we establish exponential decay estimates for the corresponding elementary symmetric functions. For characters above the threshold, we use a Gram-matrix inequality to control their number. Moreover, our character-orthogonality identity incorporates additional weights arising from character sums over effective divisors. Combining these estimates, we show that the total contribution of the nontrivial characters cannot cancel that of the trivial character when the code length is sufficiently large. This yields our square-root length bounds for sufficiently large dimensions. For smaller dimensions, we further combine higher-moment estimates with shortening arguments to obtain nontrivial length bounds beyond the range covered by the large-character method.

\subsection*{Structure of the paper}
In Section~\ref{Sect: RR_dual}, we recall the necessary background on algebraic geometry codes, establish properties of shortened codes, and prove the key obstruction Lemma~\ref{lem:obstruction}. In Section~\ref{Sect: char_sum_L-fun}, we develop the character-sum tools and establish several estimates used throughout the paper. In Section~\ref{sec:moments}, we combine these results with a character-orthogonality identity to prove our main theorems and corollaries. Finally, in Section~\ref{sec:construction}, we investigate the sharpness of our bounds through two explicit constructions.

\subsection*{Acknowledgement} 
The authors thank Yuan Ren for helpful discussions. The second author is supported by the National Natural Science Foundation of China (Grant No.~12571578).
	
\subsection*{Use of AI}
AI tools were used solely to check for mathematical gaps and grammatical errors. The authors take full responsibility for the content and accuracy of this paper.

	\section{Preliminaries}\label{Sect: RR_dual}
	In this section, we recall the Riemann--Roch and duality results used below,
	and establish certain key lemmas that will be used later. 
	
	For a divisor $H$, let $\ell(H)=\dim_{\F_q}L(H)$. Let $D$ and $G$ be as in Section~\ref{Sect: intro}. The \emph{evaluation map} is defined by 
	\begin{equation}\label{Eqn: ev}
		{\rm ev} : L(G)\to \F_q^n, \qquad f \mapsto (f(P_1), \dots, f(P_n)).
	\end{equation}
	It is immediate that we have the short exact sequence
	\begin{equation}\label{Eqn: exact}
		0\longrightarrow L(G-D)
		\longrightarrow L(G)
		\xrightarrow{\rm ev}\CL(D,G)
		\longrightarrow0.
	\end{equation}
	To simplify notations, in this paper, we always write $C=\CL(D,G)$. And denote
	$k=\dim_{\F_q}\CL(D,G)=\ell(G)-\ell(G-D)$. In particular, if $\deg G<n$, then
	$L(G-D)=\{0\}$ and the evaluation map is injective
	and $k=\ell(G)$.

	To treat MDS and NMDS codes uniformly, throughout this paper, we take $\delta\in\{0,1\}$. The codes we concern are those satisfying
	\begin{equation}\label{eq:defects}
		d(C)\ge n-k+1-\delta,
		\qquad
		d(C^\perp)\ge k+1-\delta.
	\end{equation}
	For $\delta=0$, these conditions characterize MDS codes;
	For $\delta=1$, both Singleton defects are at most one.
	Since the dual of an MDS code is MDS, these conditions
	characterize codes that are either MDS or NMDS.

	\begin{lemma}\cite[Theorem~1.6.13]{Stichtenoth}\label{Lem: Clifford}
		For all divisors $H$ with $0\leq \deg H\leq 2g-2$, one has 
		\[
		\ell(H)\leq 1+\deg(H)/2.
		\]
	\end{lemma}

	\begin{lemma}\label{lem:normalization}
		Let $g\ge2$, and let $C$ have length $n$ and dimension $k$. If $g+1\le k\le n/2$, then
		\begin{equation}\label{eq:degG}
			2g-2<\deg G=k+g-1<n.
		\end{equation}
		In particular, the evaluation map \eqref{Eqn: ev} is injective and $k=\ell(G)$.
	\end{lemma}
	
	\begin{proof}
		Suppose first that $\deg G\ge n$. Since $n\ge2k\ge2g+2$, we have $\deg G>2g-2$, and hence $\ell(G)=\deg G+1-g$ by Riemann--Roch. Using $\ell(G-D)\le\deg(G-D)+1$, we obtain
		\[
		n-k=n-\ell(G)+\ell(G-D)=g-1-\deg(G-D)+\ell(G-D)\le g.
		\]
		This contradicts $n-k\ge k\ge g+1$. Therefore $\deg G<n$, so the evaluation map is injective and $k=\ell(G)$.
		
		Suppose now that $0\le \deg G\le2g-2$. Lemma~\ref{Lem: Clifford} gives
		$
		\ell(G)\le1+\frac{\deg G}{2}\le g.
		$
		But this contradicts $k=\ell(G)\ge g+1$. Thus $\deg G>2g-2$, and Riemann--Roch yields $k=\ell(G)=\deg G+1-g$. Consequently, $\deg G=k+g-1<n$, as claimed.
	\end{proof}

	For a linear code $C\subseteq\F_q^n$, its \emph{dual code} is
	$$
	C^\perp
	=
	\left\{
	v\in\F_q^n:
	\sum_{i=1}^n c_i v_i=0
	\text{ for every }c=(c_1,\ldots,c_n)\in C
	\right\}.
	$$
	If $C$ has dimension $k$, then linear algebra gives
	$
	\dim C^\perp=n-k,
	(C^\perp)^\perp=C.
	$
	
	\begin{lemma}\cite[Proposition~2.2.10]{Stichtenoth}\label{Lem: dual_AG}
		Let $C=\CL(D,G)$. There exists an $\F_q$-rational meromorphic differential $\eta$ with a simple pole and residue $1$ at each $P_i$, and
		$
		C^\perp=\CL(D,D-G+(\eta)).
		$
		In particular, $C^\perp$ is also an algebraic geometry code.
	\end{lemma}

	Although the dimension of the dual is always determined by that of $C$, its minimum distance is generally not determined by the parameters of $C$. However, for MDS codes, the dual distance can be  determined directly.

	\begin{lemma}\cite[Chapter~11, \S\,2, Theorem~2]{MS}\label{Lem: dual_MDS}
		If $C$ is an $[n,k,n-k+1]_q$ code with $1\le k<n$, then $C^\perp$ is an $[n,n-k,k+1]_q$ code. In particular, $C$ is MDS if and only if $C^\perp$ is MDS.
	\end{lemma}
	
	\begin{remark}\label{Lem: dual_NMDS}
		By definition, a linear code $C$ is NMDS if and only if $C^\perp$ is NMDS.
	\end{remark}

	The next concept will find to be useful in our proofs. 
	\begin{definition}\label{Defn: shortened_code}
		For the divisor $D$ introduced in Section~\ref{Sect: intro}, let $E$ be an effective subdivisor of $D$. The corresponding \emph{shortened} code is defined to be 
		$$
		C_{\mathrm{short}}=\CL(D-E,G-E).
		$$
	\end{definition}
	Equivalently, consider the Riemann-Roch subspace $L(G-E)$, which consists precisely of the functions in $L(G)$ that vanish
	at every point of $\supp E$. Evaluating these functions
	at $D-E$ gives the shortened code.

	\begin{lemma}\label{lem:shortening}
		With the above notations,
		let $g\ge2$, let $\delta\in\{0,1\}$, and let
		$C$ have length
		$n$ and dimension $k$, where $g+2\le k\le n/2$.
		Assume that $C$ satisfies \eqref{eq:defects}, i.e., 
		$$
		d(C)\ge n-k+1-\delta,
		\qquad
		d(C^\perp)\ge k+1-\delta.
		$$
		Suppose $\deg E=k-m$ with $g+2\le m\le k$. Then the corresponding shortened code has length
		$s=n-k+m$ and  dimension $m$. Its distances satisfy
		$$
		d(C_{\mathrm{short}})\ge s-m+1-\delta,
		\qquad
		d(C_{\mathrm{short}}^\perp)\ge m+1-\delta.
		$$
		Moreover, $n\le2(s-m)$.
	\end{lemma}

    \begin{proof}
The length of $C_{\mathrm{short}}$ follows from its definition. Lemma~\ref{lem:normalization} applied to $C$ gives $\deg G=k+g-1<n$. Therefore the evaluation map of $C_{\mathrm{short}}$ is injective, because its kernel is $L(G-D)$ and $\deg(G-D)<0$. On the other hand, $\deg(G-E)=\deg G-(k-m)=m+g-1$. As $m\ge g+2$, we have $\deg(G-E)\ge2g+1>2g-2$, so Riemann--Roch gives $\ell(G-E)=\deg(G-E)+1-g=m$. Thus $\dim C_{\mathrm{short}}=\ell(G-E)=m$.

Shortening does not decrease the minimum distance, so
$
d(C_{\mathrm{short}})\ge d(C)\ge s-m+1-\delta.
$
The dual of the shortened code is obtained from $C^\perp$ by deleting the coordinates in $\operatorname{Supp}E$. Deleting $k-m$ coordinates lowers the minimum distance by at most $k-m$, so
$
d(C_{\mathrm{short}}^\perp)\ge d(C^\perp)-(k-m)\ge m+1-\delta.
$
Finally, $s-m=n-k\ge n/2$ because $k\le n/2$, so $n\le2(s-m)$.
\end{proof}

    Before leaving this section, we state a key lemma. In fact, later we will make use of it to prove that certain sums of characters are zero. 
	
	\begin{lemma}
		\label{lem:obstruction}
		Let $\delta\in\{0,1\}$, and let $C=\CL(D,G)$ have
		length $n$ and dimension $k$. Suppose that $\deg G<n$ and
		$
		d(C)\ge n-k+1-\delta.
		$
		Then, for every subset $B\subseteq\supp(D)$ of size $k+\delta$,
		$$
		L\left(G-\sum_{P\in B}P\right)=\{0\}.
		$$
	\end{lemma}
	
	\begin{proof}
		Suppose that there is a nonzero
		$f\in L\left(G-\sum_{P\in B}P\right)$, then $f\in L(G)$. Because 
		$\deg G<n$, we have $\operatorname{ev}(f)\neq 0$. However, since $\operatorname{ev}(f)\neq 0$ vanishes
		at every  $P\in B$, so
		$
		d(C)\le \operatorname{wt}(\operatorname{ev}(f))
		\le n-k-\delta,
		$
		contrary to the assumed minimum-distance bound.
	\end{proof}

	\section{Characters of the Jacobian: moments, sieve, and symmetric functions}\label{Sect: char_sum_L-fun}
        Our length bounds rest on estimates for character sums over the group $J=\mathcal J(\F_q)$ of rational points of the Jacobian. In this section we fix the notation for the Abel--Jacobi map, the characters of $J$ and the associated $L$-functions, and we prove the estimates used in Section~\ref{sec:moments}. The basic setup, which encodes the MDS condition as a restricted subset-sum problem in $J$, follows Chen, Hu and Ma~\cite{CHM}. The new ingredients here are the weights $W_e(\chi)$ coming from effective divisors of degree $e=g-1-\delta$, which let us treat MDS and NMDS codes uniformly, and the separate treatment of characters with large and small $|S_A(\chi)|$, which replaces the uniform bound $|S_A(\chi)|\le 2(g-1)\sqrt q$ used in~\cite{CHM}.

	Fix a rational point $P_0\in X(\F_q)$, and let $\mathcal J=\operatorname{Jac}(X)$
	be the Jacobian of $X$. Write $J=\mathcal J(\F_q)$ for the finite abelian group
	of rational points of the Jacobian. The base point $P_0$ defines the Abel--Jacobi morphism
	$
	\iota:X\longrightarrow \mathcal J$ which sends a point $
	P$ to the divisor class $[P-P_0]$. 
	We denote $h=\# J$ and $N=\#X(\F_q)$. Let $\wh J=\operatorname{Hom}(J,\mathbb C^\times)$. For $\chi\in\wh J$ and a rational divisor $R$, put
	$$
	\widetilde\chi(R)=\chi([R-(\deg R)P_0]).
	$$
	For each non-negative integer $j$, let the set $A_j$ consist of all effective rational divisors of degree $j$, including divisors containing closed points of higher degree. Denote $a_j=\#A_j$ and define 
	\begin{equation}\label{Eqn: Wj}
		W_j(\chi)=\sum_{R\in A_j}\widetilde{\chi}(R)=\sum_{R\in A_j}
		\chi([R-jP_0]).
	\end{equation}
	In particular, $W_j(1)=a_j$.

	For a subset $A\subseteq X(\F_q)$, define
	$$
	S_A(\chi)=\sum_{P\in A}\chi(\iota(P))=\sum_{P\in A}\widetilde{\chi}([P-P_0])
	$$
	In particular, $S_{X(\F_q)}(\chi)=W_1(\chi)$.

	The $L$-function associated with $\chi$ is
	$$
	L(X,\chi,t)=\prod_{P\in|X|}\left(1-\widetilde\chi(P)t^{\deg P}\right)^{-1},
	$$
	where $|X|$ denotes the set of closed points of $X$. For the trivial character, this is the zeta function $Z_X(t)$.
	Expanding the Euler product as a formal power series gives
	\begin{equation}\label{eq:character-L-series}
		L(X,\chi,t)=\sum_{j\ge0}W_j(\chi)t^j.
	\end{equation}
	
	Let $F$ be the $q$-power Frobenius endomorphism of $\mathcal J$. Pulling back the isogeny $F-1:\mathcal J\to\mathcal J$ along $\iota$ gives a geometrically connected finite \'etale abelian covering $\pi:Y\to X$ in which $P_0$ splits completely. By unramified class field theory, the natural identification
	$
	\theta:J\xrightarrow{\sim}\operatorname{Gal}(Y/X)
	$
	satisfies
	$
	\theta([P-dP_0])=\operatorname{Frob}_P
	$
	for every closed point $P$ of degree $d$, where $\operatorname{Frob}_P$ denotes arithmetic Frobenius; see \cite[Chapter~VI]{SerreClassFields}.
	
	For $\chi\in\wh J$, let $\rho_\chi=\chi\circ\theta^{-1}$ be the corresponding one-dimensional character of $\operatorname{Gal}(Y/X)$. Then
	$
	\rho_\chi(\operatorname{Frob}_P)=\widetilde\chi(P),
	$
	so $L(X,\chi,t)$ is the Artin $L$-function associated with $\rho_\chi$. Since the covering is everywhere unramified, the Artin conductor of $\rho_\chi$ is zero.

	\begin{lemma}\label{lem:character-L}
		Suppose $g\ge2$ and $\chi\in\wh J$ is nontrivial. Then
		$$
		|W_j(\chi)|\le\binom{2g-2}{j}q^{j/2}
		\qquad(0\le j\le2g-2),
		$$
		and $W_j(\chi)=0$ for $j>2g-2$. In particular,
		$
		|S_{X(\F_q)}(\chi)|\le2(g-1)\sqrt q.
		$
	\end{lemma}
	\begin{proof}
		Since $Y\to X$ is geometrically connected and everywhere unramified, the degree formula and Riemann hypothesis for Artin $L$-functions
		\cite[Theorem~A, pp.~87--88]{Bombieri} give $L(X,\chi,t)=\prod_{\nu=1}^{2g-2}(1-\alpha_\nu t)$ with $|\alpha_\nu|=\sqrt q$. 
		Comparing coefficients with
		\eqref{eq:character-L-series} gives the bounds on
		$W_j(\chi)$ and its vanishing for $j>2g-2$.
		Finally, $S_{X(\F_q)}(\chi)=W_1(\chi)$ gives the
		character-sum bound.
	\end{proof}
	
	\begin{lemma}\label{lem:moments}
		For every $A\subseteq X(\F_q)$ of size $s$ and every $0<j\le2g$, 
        \begin{equation}\label{eq:allmoments}
			\frac1h\sum_{\chi\in\wh J}|S_A(\chi)|^j\le(4g-2)^{j/2}s^{j/2}.
		\end{equation}

	\end{lemma}
	
	\begin{proof}
    By H\"older's inequality, it is sufficient to prove that for every $1\le r\le g$ we have 
    		\begin{equation}\label{eq:moments}
			\frac1h\sum_{\chi\in\wh J}|S_A(\chi)|^{2r}
			\le(4g-2)^r s^r.
		\end{equation}	
		By definition, $S_A(\chi)=\sum_{P\in A}\chi(\iota(P))=\sum_{P\in A}\widetilde{\chi}(P)$. Hence
		$
		S_A(\chi)^r=\sum_{(P_1,\ldots,P_r)\in A^r}\widetilde{\chi}\left(\sum_{i=1}^r P_i\right).
		$
		Since $\widetilde{\chi}(-R)=\overline{\widetilde{\chi}(R)}$ for every divisor $R$, we obtain
		\begin{equation}\label{Eqn: char_inner_prod}
			|S_A(\chi)|^{2r}=S_A(\chi)^r\overline{S_A(\chi)^r}=\sum_{(P_1,\ldots,P_r)\in A^r}\sum_{(Q_1,\ldots,Q_r)\in A^r}\widetilde{\chi}\left(\sum_{i=1}^r P_i-\sum_{j=1}^r Q_j\right).
		\end{equation}
		By the orthogonality relations for characters of the finite abelian group $J$, we have
		\begin{equation}\label{Eqn: orthogonality}
			\sum_{\chi\in\widehat J}\widetilde{\chi}\left(\sum_{i=1}^r P_i-\sum_{j=1}^r Q_j\right)=
			\begin{cases}
				h,&\text{if }\sum_{i=1}^r P_i\sim\sum_{j=1}^r Q_j,\\
				0,&\text{otherwise}.
			\end{cases}
		\end{equation}
		Consequently, the left-hand side of \eqref{eq:moments} equals the number of ordered pairs of tuples $((P_i)_{i=1}^r,(Q_j)_{j=1}^r)\in A^r\times A^r$ satisfying
		$
		P_1+\cdots+P_r\sim Q_1+\cdots+Q_r,
		$
		where repetitions are allowed.
		
		For a given non-negative integer $u$, we first estimate the number of ordered tuples $(P_1,\ldots,P_r)\in A^r$ for which $\ell(F)=u+1$, where $F=P_1+\cdots+P_r$. Consider the decreasing sequence of Riemann--Roch spaces
		\[
		L(K_X)\supseteq L(K_X-P_1)\supseteq L(K_X-P_1-P_2)\supseteq\cdots\supseteq L(K_X-F).
		\]
		At each step, the dimension either decreases by one or remains unchanged. We call the former a \emph{lowering step} and the latter a \emph{non-lowering step}, respectively.
		
		Recall that $\ell(K_X)=g$. Since $\deg F=r$ and $\ell(F)=u+1$, Riemann--Roch gives
		$
		\ell(K_X)-\ell(K_X-F)=r-u.
		$
		It follows that the dimension lowers at exactly $r-u$ steps, while it remains the same at the other  $u$ steps.
		
		We now estimate the number of possible choices at each step. For a lowering step, there are at most $s=|A|$ choices for $P_i$. On the other hand, for a non-lowering step, suppose that $P_1,\ldots,P_{i-1}$ have already been chosen, and put
		$
		H=K_X-P_1-\cdots-P_{i-1}.
		$
		If the $i$-th step is non-lowering, then $\ell(H-P_i)=\ell(H)$, which means that $P_i$ is a base point of the complete linear system $|H|$. Choose a nonzero $f\in L(H)$. Since $P_i$ is a base point of $|H|$, it must belong to the support of the effective divisor $\operatorname{div}(f)+H$. The degree of this divisor is
		$
		\deg H=2g-2-(i-1)=2g-1-i.
		$
		Consequently, at the $i$-th step, there are at most $2g-1-i$ choices for $P_i$ if that step is non-lowering.

		Put
		$
		c_i=2g-1-i
		$
		for 
		$
		1\le i\le r.
		$
		Suppose that the $u$ non-lowering steps occur at positions $1\le i_1<\cdots<i_u\le r$, the number of possible tuples with this prescribed pattern is at most
		$
		s^{r-u}\prod_{j=1}^u c_{i_j}.
		$
		Summing over all possible positions of the $u$ non-lowering steps, we conclude that the number $N_P(u)$ of first tuples satisfying $\ell(F)=u+1$ is bounded by
		\[
		N_P(u)\le s^{r-u}\sum_{1\le i_1<\cdots<i_u\le r}\prod_{j=1}^u c_{i_j}.
		\]
		By Maclaurin's inequality, and the fact that 
		$
		\frac1r\sum_{i=1}^r c_i=2g-1-\frac{r+1}{2},
		$
		we obtain
		\begin{equation}\label{eq:moment-first}
			N_P(u)\le\binom ru s^{r-u}\left(2g-1-\frac{r+1}{2}\right)^u.
		\end{equation}
		
		Next, fix an ordered tuple $(P_1,\ldots,P_r)\in A^r$ with $\ell(F)=u+1$, we estimate the number of ordered tuples $(Q_1,\ldots,Q_r)\in A^r$ satisfying
		$
		Q_1+\cdots+Q_r\sim F.
		$
		For such a tuple, consider the decreasing sequence
		\[
		L(F)\supseteq L(F-Q_1)\supseteq L(F-Q_1-Q_2)\supseteq\cdots\supseteq L(F-Q_1-\cdots-Q_r).
		\]
		Since $Q_1+\cdots+Q_r\sim F$, we have $F-Q_1-\cdots-Q_r\sim0$, and therefore
		$
		\ell(F-Q_1-\cdots-Q_r)=\ell(0)=1.
		$
		The initial dimension is $\ell(F)=u+1$, while the final dimension is one. Hence, the total decrease in dimension is $u$. Thus exactly $u$ steps are lowering steps, while the remaining $r-u$ steps are non-lowering steps.
		
		By similar arguments as the above, at each lowering step, there are at most $s$ choices for $Q_i$. For a non-lowering step at position $i$, the base point argument gives at most $r-i+1$ choices for $Q_i$ if that step is non-lowering.
		
		Put
		$
		d_i=r-i+1
		$
		for $1\le i\le r$. 
		Suppose that the $r-u$ non-lowering steps occur at positions $1\le j_1<\cdots<j_{r-u}\le r$. By the preceding argument, the number of second tuples with this prescribed pattern is at most
		$
		s^u\prod_{k=1}^{r-u}d_{j_k}.
		$
		Summing over all possible positions of the non-lowering steps, we obtain
		\[
		N_Q(F)\le s^u\sum_{1\le j_1<\cdots<j_{r-u}\le r}\prod_{k=1}^{r-u}d_{j_k},
		\]
		where $N_Q(F)$ denotes the number of admissible second tuples for the fixed first tuple. Since
		$
		\frac1r\sum_{i=1}^r d_i=\frac{r+1}{2},
		$
		Maclaurin's inequality gives
		\begin{equation}\label{eq:moment-second}
			N_Q(F)\le\binom ru s^u\left(\frac{r+1}{2}\right)^{r-u}.
		\end{equation}
		
		Combining \eqref{eq:moment-first} and \eqref{eq:moment-second}, we find that the number of admissible pairs of tuples for which the first tuple satisfies $\ell(F)=u+1$ is at most
		\[
		\binom ru^2s^r\left(2g-1-\frac{r+1}{2}\right)^u\left(\frac{r+1}{2}\right)^{r-u}.
		\]
		Put
		$
		\alpha=2g-1-\frac{r+1}{2},\beta=\frac{r+1}{2}.
		$
		Summing over all possible values $0\le u\le r$ and using \eqref{Eqn: char_inner_prod} and \eqref{Eqn: orthogonality}, we obtain
		\[
		\frac1h\sum_{\chi\in\wh J}|S_A(\chi)|^{2r}\le s^r\sum_{u=0}^r\binom ru^2\alpha^u\beta^{r-u}\le s^r(\sqrt{\alpha}+\sqrt{\beta})^{2r}\le s^r(2(\alpha+\beta))^r=(4g-2)^rs^r.
		\]
		This gives \eqref{eq:moments} as desired. 
	\end{proof}
	
	For the estimates below, let $g\ge\delta+2$, and put
	\begin{equation}\label{eq:aTU}
		e=g-1-\delta,\qquad T=\binom{2g-2}{e}q^{e/2},\qquad U=2(g-1)\sqrt q.
	\end{equation}
	
	\begin{lemma}\label{lem:estimates}
		Suppose $q\ge64g^2$. Then
		$
			a_e\ge\frac34q^e
		$, 
		$
			\frac{hT}{a_e}\le\frac{23}{22}\binom{2g-2}{e}q^{(g+1+\delta)/2}
		$, 
		and 
		$
			\frac{NT}{a_eU}\le 1.
		$
	\end{lemma}
	
	\begin{proof}
		We first prove the lower bound on $a_e$. For $r\ge1$, write $N_r=\#X(\mathbb F_{q^r})$. According to Weil bounds, we have 
		$
		N_r\ge M_r:=q^r+1-2gq^{r/2}>0.
		$
		Recall from the definition of Zeta function, 
		\[
		Z_X(t)=\exp\left(\sum_{r\ge1}\frac{N_rt^r}{r}\right)=\prod_{r\ge1}\sum_{j\ge0}\frac{N_r^j}{r^jj!}t^{rj}=\frac{P_X(t)}{(1-t)(1-qt)}.
		\]
		On the other hand, we define 
		\[
		Z_0(t)=\exp\left(\sum_{r\ge1}\frac{M_rt^r}{r}\right)=\prod_{r\ge1}\sum_{j\ge0}\frac{M_r^j}{r^jj!}t^{rj}=\frac{(1-\sqrt q\,t)^{2g}}{(1-t)(1-qt)}.
		\]
		Since $N_r\ge M_r$, we find that the coefficient of $Z_X(t)$ for $t^n$ is greater than or equal to that of $Z_0(t)$ for every $n$. Recall that by \eqref{eq:character-L-series},
		$
		Z_X(t)=\sum_{r\geq 0} a_rt^r
		$. 
		Hence in particular, we have $a_e$ greater than or equal to the coefficient of $t^e$ in $Z_0(t)$. 
		
		For every integer $r\ge0$, denote
		$
		S_r=\sum_{j=0}^{\min\{r,2g\}}(-1)^j\binom{2g}{j}q^{-j/2}.
		$
		Since $2g/\sqrt{q}<1/4$, the sequence $\binom{2g}{j}q^{-j/2}$ is decreasing. So the alternating-series estimate gives $S_r\ge1-2g/\sqrt q>0$ for every $r\ge0$. 
		
		Notice that 
		$
		\frac{(1-\sqrt q\,t)^{2g}}{1-qt}=\sum_{r\ge0}q^rS_rt^r.
		$
		Since $(1-t)^{-1}=\sum_{r\ge0}t^r$, it follows that
		$
		a_e\ge \sum_{r=0}^e q^rS_r.
		$
		Consequently,
		\begin{equation}\label{Eqn: 05}
			a_e\ge\sum_{r=0}^e q^rS_r\ge\left(1-\frac{2g}{\sqrt q}\right)\sum_{r=0}^e q^r=\left(1-\frac{2g}{\sqrt q}\right)\frac{q^{e+1}-1}{q-1}\ge\frac34q^e.
		\end{equation}
		This gives the claimed lower bound for $a_e$.
		
		Next, write $P_X(t)=\sum_{j=0}^{2g}p_jt^j$. By the functional equation for the $L$-polynomial \cite[Theorem~5.1.15(b)]{Stichtenoth}, we have $P_X(t)=q^gt^{2g}P_X(1/(qt))$. Since $h=P_X(1)$ by \cite[Theorem~5.1.15(c)]{Stichtenoth}, it follows that
		\begin{equation}\label{eqn: 01}
			\frac h{q^g}=P_X(1/q)=\sum_{j=0}^{2g}\frac{p_j}{q^j}.
		\end{equation}
		On the other hand, expanding $Z_X(t)=P_X(t)/((1-t)(1-qt))$ and comparing the coefficients of $t^e$, we obtain
		\begin{equation}\label{Eqn: 04}
			\frac{a_e}{q^e}=\sum_{j=0}^e\frac{p_j}{q^j}\sum_{r=0}^{e-j}q^{-r}.
		\end{equation}
		
		By the Riemann hypothesis for curves
		\begin{equation}\label{Eqn: 06}
			\frac{|p_j|}{q^j}\le\binom{2g}{j}q^{-j/2}\le\frac1{j!}\left(\frac{2g}{\sqrt q}\right)^j.
		\end{equation}
		Since $2g/\sqrt q\le1/4$, we have
		\begin{equation}\label{Eqn: 02}
			\sum_{j\ge2}\frac1{j!}\left(\frac{2g}{\sqrt q}\right)^j\le\frac{(2g/\sqrt q)^2}{2(1-2g/(3\sqrt q))}\le\frac3{88}.
		\end{equation}
		Consequently,
		$
		\sum_{j\ge1}\frac1{j!}\left(\frac{2g}{\sqrt q}\right)^j\le\frac{25}{88}.
		$

		We now estimate $hT/a_e$. By \eqref{eqn: 01} and \eqref{Eqn: 04}, and using $p_0=1$, we obtain
		$$
		\begin{aligned}
			\frac h{q^g}-\frac{a_e}{q^e}
			&\le\sum_{j=e+1}^{2g}\frac{|p_j|}{q^j}-\frac1q+\frac1{q-1}\sum_{j=1}^e\frac{|p_j|}{q^j}\\
			&\le\sum_{j\ge2}\frac1{j!}\left(\frac{2g}{\sqrt q}\right)^j-\frac1q+\frac1{q-1}\sum_{j\ge1}\frac1{j!}\left(\frac{2g}{\sqrt q}\right)^j\\
			&\le\frac3{88}-\frac1q+\frac{25}{88(q-1)}
			\le\frac3{88}.
		\end{aligned}
		$$
		where the last inequality follows from $q\ge256$. Since $a_e/q^e\ge3/4$ by \eqref{Eqn: 05}, it follows that
		$$
		\frac h{q^g}\le\frac{a_e}{q^e}+\frac3{88}\le\frac{23}{22}\frac{a_e}{q^e}.
		$$
		Multiplying by $T=\binom{2g-2}{e}q^{e/2}$ and using $e=g-1-\delta$ gives bound for $\frac{hT}{a_e}$.

        It remains to prove the bound for $\frac{NT}{a_eU}$. If $e=1$, then $a_e=N$ and $T=U$, so $\frac{NT}{a_eU}=1$. Assume $e\ge2$. Since $\sqrt q\ge8g$, the Hasse--Weil bound gives $N\le q+1+2g\sqrt q\le q+1+\frac q4\le\frac{63}{50}q$, where the last step uses $q\ge256$. Together with $a_e\ge\frac34q^e$ this yields $\frac N{a_e}\le\frac{63}{50}\cdot\frac43\,q^{1-e}<1.7\,q^{1-e}$. Using $\binom{2g-2}{e}\le(2g-2)^e/e!$ and $\sqrt q\ge8g$, we obtain
\[
\frac{NT}{a_eU}\le1.7\,q^{1-e}\cdot\frac{\binom{2g-2}{e}q^{e/2}}{2(g-1)\sqrt q}\le\frac{1.7}{e!}\Bigl(\frac{2g-2}{\sqrt q}\Bigr)^{e-1}\le\frac{1.7}{e!\,4^{e-1}}\le\frac{1.7}{8}<1.
\]
Thus $NT/(a_eU)\le1$ for every $e\ge1$. This completes the proof.
	\end{proof}

	\begin{lemma}\label{lem:sieve}
		Let $A\subseteq X(\F_q)$ have size $s$, and let $\Gamma\subseteq\wh J$ be a
		nonempty subset of size $R$.
		\begin{enumerate}[label=\textup{(\alph*)},leftmargin=2.1em]
			\item $\displaystyle\sum_{\chi\in\Gamma}|S_A(\chi)|^2\le s\bigl(N+(R-1)U\bigr).$
			\item If $N\ge U$ and $t^2>sU$, then the number of $\chi\in\wh J$ with
			$|S_A(\chi)|\ge t$ is at most $\dfrac{s(N-U)}{t^2-sU}$.
		\end{enumerate}
	\end{lemma}
	
	\begin{proof}
		(a) For $\chi\in\Gamma$ let $v_\chi\in\mathbb C^{X(\F_q)}$ be the vector with
		coordinates $v_\chi(P)=\chi(\iota(P))$, and let $x=\mathbf 1_A$. Then
		$\|x\|^2=s$ and $|\langle x,v_\chi\rangle|=|S_A(\chi)|$. Moreover
		$\|v_\chi\|^2=N$, and for $\chi\ne\psi$ Lemma~\ref{lem:character-L}, applied to
		the nontrivial character $\chi\psi^{-1}$, gives
		$|\langle v_\chi,v_\psi\rangle|=|S_{X(\F_q)}(\chi\psi^{-1})|\le U$. Hence every row sum of the matrix $\mathcal G=(\langle v_\chi,v_\psi\rangle)_{\chi,\psi\in\Gamma}$
		is at most $N+(R-1)U$. Let $\lambda$ be an eigenvalue of the Hermitian matrix
		$\mathcal G$ with eigenvector $u\neq0$, and choose $i$ with $|u_i|=\max_j|u_j|$. Then
		$
		|\lambda||u_i|=\Bigl|\sum_j\mathcal G_{ij}u_j\Bigr|\le|u_i|\sum_j|\mathcal G_{ij}|,
		$
		so $\|\mathcal G\|=\max|\lambda|\le N+(R-1)U$. Let $V$ be the matrix with columns
		$v_\chi$, so $\mathcal G=V^*V$ up to transposition.
		So we have 
		\[
		\sum_{\chi\in\Gamma}|\langle x,v_\chi\rangle|^2=\|V^*x\|^2\le\|VV^*\|\,\|x\|^2
		=\|\mathcal G\|\,\|x\|^2\le s\bigl(N+(R-1)U\bigr).
		\]

		(b) Let $\Gamma_t=\{\chi:|S_A(\chi)|\ge t\}$ and $R_t=\#\Gamma_t$. If $R_t=0$ the
		bound is trivial, since $N\ge U$. Otherwise (a) gives
		$R_tt^2\le s(N+(R_t-1)U)$, i.e.\ $R_t(t^2-sU)\le s(N-U)$.
	\end{proof}

	For $0\le r\le s=|A|$, let
	\[
	E_r(\chi)=\sum_{\substack{B\subseteq A\\|B|=r}}
	\chi\left(\sum_{P\in B}\iota(P)\right).
	\]
    Then one can check that $E_j$ is the coefficient of $z^j$ in
		$
		F(z)=\prod_{P\in A}(1+z\chi(\iota(P))).
		$
	
	\begin{lemma}\label{lem:coefficients}
		Let $A\subseteq X(\F_q)$ have cardinality $s$. 
		Let $1\le r\le s-1$ and $\lambda=\min\{r,s-r\}$. If $0\le\theta\le1$ and $|S_A(\chi)|\le(1-\theta)s$, then
		\begin{equation}\label{eq:Cauchy-general}
			\frac{|E_r(\chi)|}{\binom sr}\le(s+1)e^{-\theta\lambda/2}.
		\end{equation}
		In particular, if $|S_A(\chi)|\le s/2$, then
		\begin{equation}\label{eq:Cauchy}
			\frac{|E_r(\chi)|}{\binom sr}\le(s+1)e^{-\lambda/4}.
		\end{equation}
	\end{lemma}
	
	\begin{proof}
		Put $p=r/s\in(0,1)$ and $\rho=p/(1-p)=r/(s-r)$. We write $\Re(z)$ for the real part of a complex number $z$. Since $|\chi(\iota(P))|=1$ for every $P\in A$, the arithmetic--geometric mean inequality gives, for every $\varphi\in\mathbb R$,
		\begin{align*}
			|F(\rho e^{i\varphi})|^2 &=\prod_{P\in A}|1+\rho e^{i\varphi}\chi(\iota(P))|^2\\
			&\le\left(\frac1s\sum_{P\in A}|1+\rho e^{i\varphi}\chi(\iota(P))|^2\right)^s\\
			&=\left(1+\rho^2+\frac{2\rho}{s}\Re(e^{i\varphi}S_A(\chi))\right)^s\\
			&\le\left((1+\rho)^2-2\rho\theta\right)^s.
		\end{align*}
		Using $\rho/(1+\rho)^2=p(1-p)$ and $1-u\le e^{-u}$, we obtain
		$
		|F(\rho e^{i\varphi})|\le(1+\rho)^s e^{-\theta sp(1-p)}.
		$
		Cauchy's coefficient estimate therefore yields
		$$
		|E_r(\chi)|\le\frac{e^{-\theta sp(1-p)}}{p^r(1-p)^{s-r}}.
		$$
		To normalize this bound, take
		$
		b_j=\binom{s}{j}p^j(1-p)^{s-j}$ for $0\le j\le s.
		$
		The term $b_r$ is a largest term of this sequence. As $\sum_{j=0}^s b_j=1$, we have $b_r\ge1/(s+1)$. Consequently,
		$$
		\frac{|E_r(\chi)|}{\binom{s}{r}}\le(s+1)\exp\left(-\theta\frac{r(s-r)}s\right).
		$$
		Since $r(s-r)/s\ge\min\{r,s-r\}/2=\lambda/2$, this proves \eqref{eq:Cauchy-general}. Taking $\theta=1/2$ gives \eqref{eq:Cauchy}.
	\end{proof}

	\section{Proofs of the main theorems}\label{sec:moments}

We begin with the idea common to all proofs. Let $C'=\CL(D',G')$ be either $C$ itself or a shortened code of $C=\CL(D,G)$ (Definition~\ref{Defn: shortened_code}). Put $A=\supp(D')$, $s=\#A$ and $e=g-1-\delta$, and suppose that $C'$ has dimension $m-\delta$, with $\deg G'=m+e<s$ and minimum distance at least $s-m+1$. By Lemma~\ref{lem:obstruction}, there is no subset $B\subseteq A$ with $\#B=m$ and no effective rational divisor $H$ of degree $e$ such that $G'\sim\sum_{P\in B}P+H$. Character orthogonality on $J$ turns this non-existence into the identity
\begin{equation}\label{eq:common-character-identity}
0=\frac1h\sum_{\chi\in\wh J}\chi\bigl(-[G'-(m+e)P_0]\bigr)E_m(\chi)W_e(\chi),
\end{equation}
where $E_m(\chi)$ and $W_e(\chi)$ are the character sums over $m$-subsets of $A$ and over effective divisors of degree $e$, respectively. The trivial character contributes the positive quantity $\binom sm a_e/h$, so it must be cancelled by the nontrivial characters. For $\chi\ne1$ we make use of the estimation in the previous section, which helps us to conclude: if $s$ is large these two contributions cannot cancel the trivial one. This gives rise to the bounds for $s$ and hence $n$.

We first apply this strategy in three ranges of the dimension $k$ relative to the genus $g$. In Section~\ref{Sect: finish_proof}, we combine the resulting estimates to prove Theorems~\ref{thm:length} and~\ref{thm:nmds-length}. Section~\ref{Cor: proof_Cor} derives the asymptotic corollaries from these bounds.  Finally, Section~\ref{Sect: pf_prop:nmds-genus-two} proves the sharp upper bound for genus-two NMDS codes.

We retain the notation of Sections~\ref{Sect: RR_dual} and~\ref{Sect: char_sum_L-fun}, summarized in Table~\ref{tab:notation}. Throughout, $\delta\in\{0,1\}$ and $g\ge\delta+2$, so that $e\ge1$. Let $E$ be an effective subdivisor of $D$. Set $D'=D-E$ and $G'=G-E$. The corresponding shortened code is $\CL(D',G')$. We denote its length by $s$ and its dimension by $m-\delta$.

\begin{table}[htbp]
\centering
\renewcommand{\arraystretch}{1.25}
\begin{tabular}{@{}l p{0.72\textwidth}@{}}
\hline
Symbol & Meaning \\
\hline
$\delta$ & $0$ for MDS codes, $1$ for NMDS codes; the code satisfies \eqref{eq:defects} \\
$e=g-1-\delta$ & degree of the effective divisors $H$ in the obstruction of Lemma~\ref{lem:obstruction} \\
$N=\#X(\F_q)$ & number of rational points of $X$ \\
$J=\mathcal J(\F_q)$, $h=\#J$ & group of rational points of the Jacobian, and its order \\
$\wh J$, $\iota$ & character group of $J$, and the Abel--Jacobi map $P\mapsto[P-P_0]$ \\
$a_j$ & number of effective rational divisors of degree $j$ \\
$A$, $s$ & $A=\supp(D')$ and $s=\#A$, the length of the (shortened) code \\
$S_A(\chi)$ & $\sum_{P\in A}\chi(\iota(P))$ \\
$W_j(\chi)$ & $\sum_{R\in A_j}\chi([R-jP_0])$, a sum over effective divisors of degree $j$ \\
$E_r(\chi)$ & $\sum_{B\subseteq A,\,\#B=r}\chi\bigl(\sum_{P\in B}\iota(P)\bigr)$, the $r$-th elementary symmetric function of the values $\chi(\iota(P))$, $P\in A$ \\
$T=\binom{2g-2}{e}q^{e/2}$ & bound for $|W_e(\chi)|$, $\chi\ne1$ (Lemma~\ref{lem:character-L}) \\
$U=2(g-1)\sqrt q$ & bound for $|S_{X(\F_q)}(\chi)|$, $\chi\ne1$ (Lemma~\ref{lem:character-L}) \\
\hline
\end{tabular}
\caption{Notation used in Section~\ref{sec:moments}.}
\label{tab:notation}
\end{table}

	\begin{proposition}\label{prop:count}
		Let $\delta\in\{0,1\}$, and assume that $g\ge\delta+2$ and $g+\delta+2\le m\le2g$. Suppose that $\CL(D',G')$ has minimum distance at least $s-m+1$ and that $\deg G'=m+g-1-\delta<s$. If $a_e>0$, then
		\begin{equation}\label{eq:countbound}
			s\le\binom{m}{2}+(4g-2)\left[\left(1+\frac{\binom{m}{2}}{4g-2}\right)\frac{hT}{a_e}\right]^{2/m}.
		\end{equation}
	\end{proposition}

	\begin{proof}
		Let $A=\supp(D')$. We count pairs consisting of an ordered tuple $(Q_1,\ldots,Q_m)\in A^m$, with repetitions allowed, and an effective $\F_q$-rational divisor $R$ of degree $e$ satisfying
		\begin{equation}\label{eq:target}
			\iota(Q_1)+\cdots+\iota(Q_m)+[R-eP_0]\sim [G'-(m+e)P_0]
		\end{equation}
		By character orthogonality, the total number $Z$ of such pairs is
		$$
		Z=\frac1h\sum_{\chi\in\wh J}\chi(-[G'-(m+e)P_0])S_A(\chi)^mW_e(\chi).
		$$
		The trivial character contributes $s^ma_e/h$. For the nontrivial characters, we have 
		\begin{align*}
			\left|\frac1h\sum_{\substack{\chi\in\wh J\\\chi\ne1}}\chi(-[G'-(m+e)P_0])S_A(\chi)^mW_e(\chi)\right| &\le\frac{T}{h}\sum_{\substack{\chi\in\wh J\\\chi\ne1}}|S_A(\chi)|^m \qquad \text{by Lemma~\ref{lem:character-L}}\\
			&\le T(4g-2)^{m/2}s^{m/2} \qquad \text{by \eqref{eq:allmoments}}
		\end{align*}
		Therefore,
		\begin{equation}\label{eq:Z}
			Z\ge\frac{s^ma_e}{h}-T(4g-2)^{m/2}s^{m/2}.
		\end{equation}

		For a pair of distinct positions $i,j$, let $Z_{ij}$ count solutions of \eqref{eq:target} with $Q_i=Q_j$. Then 
		\[
		Z_{ij}=\frac1h\sum_{\chi\in\wh J}\chi(-[G'-(m+e)P_0])S_A(\chi^2)S_A(\chi)^{m-2}W_e(\chi).
		\]
		By similar arguments as the above, the trivial character contributes $s^{m-1}a_e/h$. Since $|S_A(\chi^2)|\le s$ and $m-2\ge2$, the moment bound of order $m-2$ gives
		\begin{equation}\label{eq:collision}
			Z_{ij}\le\frac{s^{m-1}a_e}{h}+T(4g-2)^{(m-2)/2}s^{m/2}.
		\end{equation}
		
		Apply Lemma~\ref{lem:obstruction} by taking $n=s$ and $k=m-\delta$, then we see that \eqref{eq:target} cannot have a solution with all $Q_i$ distinct. In fact, suppose all $Q_i$ are distinct, and there exists an effective $R$ satisfying \eqref{eq:target}, then 
		$
		\ell(G'-\sum_{i=1}^m Q_i)=\ell(R)\neq 0
		$
		which is impossible. As a consequence, every solution to \eqref{eq:target} must have at least two $Q_i$ equal, and hence $Z\le\sum_{i<j}Z_{ij}$. Combining \eqref{eq:Z} and \eqref{eq:collision} with $Z\le\sum_{i<j}Z_{ij}$, and dividing by $(a_e/h)s^{m/2}$, gives
		\begin{equation}\label{eq:before-root}
			s^{m/2-1}\left(s-\binom{m}{2}\right)\le(4g-2)^{m/2}\left(1+\frac{\binom{m}{2}}{4g-2}\right)\frac{hT}{a_e}.
		\end{equation}
		If $s\le\binom{m}{2}$, the bound \eqref{eq:countbound} is immediate. Otherwise, the left-hand side of \eqref{eq:before-root} is at least $\left(s-\binom{m}{2}\right)^{m/2}$. Raising both sides to the power $2/m$ proves \eqref{eq:countbound}.
	\end{proof}
	
\begin{proposition}\label{prop:uniform}
Let $\delta\in\{0,1\}$ and $g\ge\delta+2$, and let $C=\CL(D,G)$ have length $n$ and dimension $k$ with $g+2\le k\le n/2$ and $k+\delta\le2g$. If $C$ satisfies \eqref{eq:defects}, then
\begin{equation}\label{eq:uniform-sharp}
n\le65g\,q^{(g+1+\delta)/(k+\delta)}.
\end{equation}
\end{proposition}

\begin{proof}
Put $m=k+\delta$. Since $k\ge g+2$ and $k+\delta\le2g$, we have $g+\delta+2\le m\le2g$. If $q<64g^2$, then $\sqrt q<8g$, and the Hasse--Weil bound gives $n\le N\le q+1+2g\sqrt{q}\le 65g\,q^{(g+1+\delta)/m}$. Hence we may assume $q\ge64g^2$.

We apply Proposition~\ref{prop:count} to $C$ itself, i.e.\ with $E=0$, $s=n$, $D'=D$, $G'=G$. Then $\deg G=k+g-1=m+g-1-\delta<n$ by Lemma~\ref{lem:normalization}, $d(C)\ge n-k+1-\delta=n-m+1$ by \eqref{eq:defects}, and $a_e>0$ by Lemma~\ref{lem:estimates}, since $q\ge64g^2$. Therefore
\[
n\le\binom m2+(4g-2)\left[\Bigl(1+\frac{\binom m2}{4g-2}\Bigr)\frac{hT}{a_e}\right]^{2/m}.
\]

We estimate the expression in the square brackets. Since $m\le2g$, we have $\binom m2\le g(2g-1)\le2g^2$ and $1+\frac{\binom m2}{4g-2}\le1+\frac{g(2g-1)}{2(2g-1)}=\frac{g+2}2$. By Lemma~\ref{lem:estimates}, $\frac{hT}{a_e}\le\frac{23}{22}\binom{2g-2}{g-1-\delta}q^{(g+1+\delta)/2}$. From $\binom{2r}{r}/4^r=\prod_{j=1}^r(1-\frac1{2j})\le\exp\bigl(-\tfrac12\log(r+1)\bigr)=(r+1)^{-1/2}$ we get $\binom{2g-2}{g-1-\delta}\le\binom{2g-2}{g-1}\le4^{g-1}/\sqrt g$. Combining all these, 
\[
\Bigl(1+\frac{\binom m2}{4g-2}\Bigr)\frac{hT}{a_e}\le\frac{23}{22}\cdot\frac{g+2}{2}\binom{2g-2}{g-1-\delta}\,q^{(g+1+\delta)/2}\le4^{g+2}\cdot\frac{23(g+2)}{2816\sqrt g}\,q^{(g+1+\delta)/2}.
\]
We claim that
$
\frac{23(g+2)}{2816\sqrt g}\le\Bigl(1+\frac1{256}\Bigr)^{(g+2)/2}.
$
For $2\le g\le17$ the left-hand side is less than $1$. For $g\ge18$ we have $\sqrt g>4.24$, so the left-hand side is less than $1+\frac{g+2}{512}\le(1+\frac1{256})^{(g+2)/2}$. This proves the claim. Raising the preceding display to the power $2/m$ and using $m\ge g+2$, we obtain
\[
\begin{aligned}
    \left[\Bigl(1+\frac{\binom m2}{4g-2}\Bigr)\frac{hT}{a_e}\right]^{2/m}& \le\Bigl(4^{g+2}\bigl(1+\tfrac1{256}\bigr)^{(g+2)/2}\Bigr)^{2/m}q^{(g+1+\delta)/m}\\
    &=16^{(g+2)/m}\bigl(1+\tfrac1{256}\bigr)^{(g+2)/m}q^{(g+1+\delta)/m}\le\frac{257}{16}\,q^{(g+1+\delta)/m}.
\end{aligned}
\]

Substituting into the bound for $n$ gives $n\le2g^2+\frac{257}{16}(4g-2)\,q^{(g+1+\delta)/m}$. Finally, $q^{(g+1+\delta)/m}\ge\sqrt q\ge8g$ gives $2g^2\le\frac g4\,q^{(g+1+\delta)/m}$, and hence
\[
n\le\Bigl(\frac g4+\frac{257}{16}(4g-2)\Bigr)q^{(g+1+\delta)/m}=\Bigl(\frac{129}{2}g-\frac{257}{8}\Bigr)q^{(g+1+\delta)/m}<65g\,q^{(g+1+\delta)/m}.\qedhere
\]
\end{proof}

\begin{proposition}\label{prop:high}
Let $\delta\in\{0,1\}$, $g\ge\delta+2$, and let $C=\CL(D,G)$ have length $n$ and
dimension $k$ with $g+2\le k\le n/2$ and $k+\delta\ge2g$. If $C$ satisfies
\eqref{eq:defects}, then
\[
n\le 33g\,q^{(g+1+\delta)/(2g)}.
\]
\end{proposition}

\begin{proof}
If $q<64g^2$, then $N\le q+1+2g\sqrt q<11g\sqrt q$, and the claim follows from
the Hasse--Weil bound since $q^{(g+1+\delta)/(2g)}\ge\sqrt q$. Therefore we assume $q\ge64g^2$.

Put $m=2g$, then we have $g+\delta+2\le m\le k+\delta$. Choose an effective subdivisor $E\le D$ of degree
$k+\delta-m$ and let $\CL(D-E,G-E)$ be the shortened code. By
Lemma~\ref{lem:shortening}, it has length $s=n-k+m-\delta$, dimension $m-\delta$,
minimum distance at least $s-m+1$, and $n\le2(s-m+\delta)$. Also
$\deg(G-E)=m+g-1-\delta<s$, so Proposition~\ref{prop:count} applies. Combined with
Lemma~\ref{lem:estimates} and $1+\binom m2/(4g-2)=(g+2)/2$, it gives
\[
n\le 4g^2-6g+2\delta+(8g-4)\Bigl(\tfrac{23(g+2)}{44}\tbinom{2g-2}{g-1-\delta}\Bigr)^{1/g}
q^{(g+1+\delta)/(2g)}.
\]
Using $\binom{2g-2}{g-1-\delta}\le\binom{2g-2}{g-1}\le 4^{g-1}/\sqrt g$,
\[
\tfrac{23(g+2)}{44}\tbinom{2g-2}{g-1-\delta}\le 4^g\,\tfrac{23(g+2)}{176\sqrt g}
\le 4^g\Bigl(1+\tfrac g{64}\Bigr)\le\Bigl(\tfrac{65}{16}\Bigr)^g,
\]
where the middle step uses $\sqrt g\le g/9+9/4$ and $2/\sqrt g\le\sqrt2$, and the
last is Bernoulli's inequality. Since $q^{(g+1+\delta)/(2g)}\ge\sqrt q\ge8g$,
\[
n\le\Bigl(\tfrac{4g^2-6g+2\delta}{8g}+(8g-4)\tfrac{65}{16}\Bigr)q^{(g+1+\delta)/(2g)}
<33g\,q^{(g+1+\delta)/(2g)}.\qedhere
\]
\end{proof}

\begin{proposition}\label{prop:root}
Let $\delta\in\{0,1\}$ and $g\ge\delta+2$, $g+2\le k\le n/2$. Suppose that $C$ satisfies \eqref{eq:defects}. If
\begin{equation}\label{eq:Kdelta}
k\ge L_{g,\delta}(q):=\delta+2(g+\delta+3)\log q+4\log\left(16\binom{2g-2}{g-1-\delta}\right),
\end{equation}
then
\[
n<20(g-1)\sqrt q.
\]
\end{proposition}

\begin{proof}
If $q<64g^2$, the desired bound follows directly from the Hasse--Weil bound. We may therefore assume that $q\ge64g^2$. By Lemma~\ref{lem:normalization}, we have $\deg G=k+g-1<n$. Put $A=\supp(D)$ and $e=g-1-\delta$. Since $\deg G=k+\delta+e$, Lemma~\ref{lem:obstruction} excludes every pair consisting of a $(k+\delta)$-element subset $B\subseteq A$ and an effective rational divisor $H$ of degree $e$ satisfying
$
G\sim\sum_{P\in B}P+H.
$
Character orthogonality therefore gives
\begin{equation}\label{eq:root-count}
0=\frac1h\sum_{\chi\in\wh J}\chi\bigl(-[G-(k+\delta+e)P_0]\bigr)E_{k+\delta}(\chi)W_e(\chi).
\end{equation}
By Lemma~\ref{lem:estimates}, $a_e>0$, so the trivial character contributes the positive quantity
$
\binom n{k+\delta}a_e/h.
$

Suppose, for a contradiction, that $n\ge10U$. Call a character \emph{large} if $|S_A(\chi)|\ge n/2$, and let $R$ be the number of large characters, including the trivial character. Since $n\ge10U>4U$ and $N\ge U$, Lemma~\ref{lem:sieve}\textup{(b)} gives
\[
R\le\frac{n(N-U)}{n^2/4-nU}=\frac{4(N-U)}{n-4U}\le\frac{4N}{n-4U}.
\]
For every nontrivial character, we have $|E_{k+\delta}(\chi)|\le\binom n{k+\delta}$ and $|W_e(\chi)|\le T$ by Lemma~\ref{lem:character-L}. Hence the total absolute contribution of the nontrivial large characters to \eqref{eq:root-count}, divided by that of the trivial character, is at most
\begin{equation}\label{eq:large-part}
\frac{RT}{a_e}\le\frac{4NT}{(n-4U)a_e}\le\frac{4U}{n-4U}\le\frac23,
\end{equation}
where we used Lemma~\ref{lem:estimates} and $n\ge10U$.

Next, since $g+2\le k\le n/2$, we have
$
\lambda:=\min\{k+\delta,n-k-\delta\}\ge k-\delta.
$
Every character that is not large is nontrivial and satisfies $|S_A(\chi)|<n/2$. Lemma~\ref{lem:coefficients} therefore gives
\[
\frac{|E_{k+\delta}(\chi)|}{\binom n{k+\delta}}\le(n+1)e^{-\lambda/4}\le(n+1)e^{-(k-\delta)/4}.
\]
There are at most $h$ such characters. Since $q\ge64g^2$, the Hasse--Weil bound gives $n+1\le N+1\le2q$. Using Lemma~\ref{lem:estimates}, we find that their total absolute contribution to \eqref{eq:root-count}, divided by that of the trivial character, is at most
\begin{equation}\label{eq:small-part}
(n+1)\frac{hT}{a_e}e^{-(k-\delta)/4}\le4\binom{2g-2}{g-1-\delta}q^{(g+\delta+3)/2}e^{-(k-\delta)/4}\le\frac14,
\end{equation}
where the last inequality follows from \eqref{eq:Kdelta}.

Combining \eqref{eq:large-part} and \eqref{eq:small-part}, the total absolute contribution of the nontrivial characters is strictly smaller than the contribution of the trivial character. This contradicts \eqref{eq:root-count}. Therefore
$
n<20(g-1)\sqrt q
$
as desired.
\end{proof}

	\subsection{Proof of Theorems~\ref{thm:length} and~\ref{thm:nmds-length}}\label{Sect: finish_proof}

\begin{proof}[Proof of Theorem~\ref{thm:length}]
Since $C$ is MDS, condition \eqref{eq:defects} holds with $\delta=0$.
\begin{enumerate}
\item If $g+2\le k<2g$, the asserted bound follows from Proposition~\ref{prop:uniform} with $\delta=0$ and $m=k$.
\item If $2g\le k<L_g(q)$, the asserted bound follows from Proposition~\ref{prop:high} with $\delta=0$.
\item If $k\ge L_g(q)$, the asserted bound follows from Proposition~\ref{prop:root} with $\delta=0$.
\end{enumerate}
Together, these cases establish \eqref{eq:mds-length-bounds}.
\end{proof}

\begin{proof}[Proof of Theorem~\ref{thm:nmds-length}]
Since $C$ is NMDS, condition \eqref{eq:defects} holds with $\delta=1$. As $g\ge3$, all the relevant propositions apply.
\begin{enumerate}
\item If $g+2\le k<2g-1$, the asserted bound follows from Proposition~\ref{prop:uniform} with $\delta=1$ and $m=k+1$.
\item If $2g-1\le k<L'_g(q)$, the asserted bound follows from Proposition~\ref{prop:high} with $\delta=1$.
\item If $k\ge L'_g(q)$, the asserted bound follows from Proposition~\ref{prop:root} with $\delta=1$.
\end{enumerate}
Together, these cases establish \eqref{eq:nmds-length-bounds}.
\end{proof}

\subsection{Proofs of Corollaries~\ref{thm:mds} and~\ref{cor:nmds}}\label{Cor: proof_Cor}

We first recall several results from finite geometry concerning MDS codes. Let $C$ be a nontrivial $[n,k,n-k+1]_q$ MDS code, and let $M\in\F_q^{k\times n}$ be a generator matrix of $C$. The columns of $M$ determine $n$ distinct points in $\operatorname{PG}(k-1,q)$ forming an $n$-arc. Indeed, the MDS property is equivalent to the condition that every $k$ columns of $M$ are linearly independent, or equivalently, that no $k$ of the corresponding projective points lie in a hyperplane.

\begin{lemma}\label{lem:arcs}
Let $C$ be a nontrivial $[n,k,n-k+1]_q$ MDS code. Then the following bounds hold.
\begin{enumerate}[label=\textup{(\roman*)},leftmargin=2.1em]
\item If $k=2$, then $n\le q+1$. If $k=3$, then $n\le q+1$ for odd $q$ and $n\le q+2$ for even $q$.
\item If $k\in\{4,5\}$ and $q\ge5$, then $n\le q+1$.
\item If $q$ is odd, $q\ge529$, $q\notin\{3^6,5^5\}$, and $6\le k\le\frac12\sqrt q+2$, then $n\le q+1$.
\item If $q$ is even and $6\le k<\frac12\sqrt q+\frac{15}{4}$, then $n\le q+1$.
\end{enumerate}
In particular, if $q>3125$ and $4\le k<\sqrt q/2$, then $n\le q+1$.
\end{lemma}

\begin{proof}
If $k=2$, the corresponding arc lies in $\operatorname{PG}(1,q)$, which contains exactly $q+1$ points. Suppose that $k=3$. The corresponding arc lies in $\operatorname{PG}(2,q)$, with no three points collinear. Fixing a point of the arc, the remaining $n-1$ points lie on distinct lines through it. Since there are $q+1$ such lines, we obtain $n\le q+2$. If equality holds, every line through a point of the arc contains exactly one other point of the arc. Thus every line meets the arc in either zero or two points. Counting the points of the arc on the lines through any point outside the arc shows that $n$ is even. Since $n=q+2$, it follows that $q$ is even. This proves \textup{(i)}.

Part~\textup{(ii)} follows from \cite[Theorems~61 and~63]{BL}. Part~\textup{(iii)} is due to Hirschfeld and Korchm\'aros \cite{HK}, while part~\textup{(iv)} follows from Storme and Thas \cite{ST}. These results and their precise ranges are summarized in \cite[Section~19]{BL}.

Finally, suppose that $q>3125$ and $4\le k<\sqrt q/2$. If $k\in\{4,5\}$, the conclusion follows from \textup{(ii)}. If $k\ge6$ and $q$ is odd, then $q\ge529$ and $q\notin\{3^6,5^5\}$, so \textup{(iii)} applies. If $k\ge6$ and $q$ is even, then \textup{(iv)} applies. Hence $n\le q+1$ in all cases.
\end{proof}

\begin{proof}[Proof of Corollary~\ref{thm:mds}]
Fix $g\ge2$ and $\varepsilon>0$. Since $L_g(q)=O_g(\log q)$, we may assume that $q$ is sufficiently large such that 
$
L_g(q)<\frac{\sqrt q}{2}$ and $q^\varepsilon\ge20(g-1).
$

If $k<L_g(q)$, then $k<\sqrt q/2$, and Lemma~\ref{lem:arcs} gives $n\le q+1$, with the stated exception when $q$ is even and $k=3$. If $k\ge L_g(q)$, Theorem~\ref{thm:length} gives
$
n<20(g-1)\sqrt q\le q^{1/2+\varepsilon}.
$
This proves the asserted bounds.

Taking $\varepsilon=1/2$, we obtain the MDS conjecture for codes of genus $g\ge2$ over all sufficiently large finite fields. For genus zero, the bound follows from $n\le\#X(\F_q)=q+1$. For genus one, Han and Ren's result \cite[Theorem~1.2]{HR} gives
$
n<q+1
$
for all sufficiently large $q$. Thus the MDS conjecture holds for algebraic geometry codes of every fixed genus over all sufficiently large finite fields.
\end{proof}

\begin{proof}[Proof of Corollary~\ref{cor:nmds}]
Fix $g\ge3$ and $\varepsilon>0$, and take $q$ sufficiently large that $q\ge4g^2$ and $q^\varepsilon\ge20(g-1)$. By the Hasse--Weil bound,
$
n\le\#X(\F_q)\le q+1+2g\sqrt q\le2q+1.
$
In particular, the asserted bound holds when $k<L'_g(q)$. If $k\ge L'_g(q)$, Theorem~\ref{thm:nmds-length} gives
$
n<20(g-1)\sqrt q\le q^{1/2+\varepsilon}.
$
This establishes both bounds. 
\end{proof}

\subsection{Proof of Proposition~\ref{prop:nmds-genus-two}}\label{Sect: pf_prop:nmds-genus-two}

	\begin{proof}[Proof of the upper bound in Proposition~\ref{prop:nmds-genus-two}]
    We prove this result by contradiction. Suppose that $n>(N+U)/2$. Let $A=\supp(D)$. By Lemmas~\ref{lem:normalization} and~\ref{lem:obstruction}, $\deg G=k+1<n$ and no sum of $k+1$ distinct points of $A$ is linearly equivalent to $G$. Character orthogonality gives
\[
0=\frac1h\sum_{\chi\in\wh J}\widetilde\chi(-G)E_{k+1}(\chi).
\]
The trivial character contributes $\binom n{k+1}/h$. Thus
\begin{equation}\label{eq:genus2-count}
1\le\sum_{\chi\ne1}\frac{|E_{k+1}(\chi)|}{\binom n{k+1}}.
\end{equation}

Call a character large if $|S_A(\chi)|\ge5n/6$. There are at most two large characters. Indeed, if there were at least three, Lemma~\ref{lem:sieve}\textup{(a)} would give $25n/12\le N+2U$. But $q\ge2500$ implies $N\ge q+1-4\sqrt q>23U$, and hence $25n/12>25(N+U)/24>N+2U$, a contradiction. The trivial character is large, and complex conjugation preserves largeness. Therefore any nontrivial large character must be quadratic.

Suppose that there is a nontrivial large character $\chi$. As $\chi$ is quadratic, its values on $A$ are $\pm1$. Both signs must occur, since Lemma~\ref{lem:character-L} gives
$
|S_A(\chi)|\le U+N-n<n.
$

Choose $P_+,P_-\in A$ such that $\chi(\iota(P_+))=1$ and $\chi(\iota(P_-))=-1$. For every $k$-element subset $B\subseteq A\setminus\{P_+,P_-\}$, the two $(k+1)$-element subsets $B\cup\{P_+\}$ and $B\cup\{P_-\}$ have opposite character values. That says, among the $\binom n{k+1}$ subsets of $A$ of cardinality $k+1$, at least $\binom{n-2}{k}$ have character value $1$, and at least $\binom{n-2}{k}$ have character value $-1$. Therefore,
\[
\frac{|E_{k+1}(\chi)|}{\binom n{k+1}}\le1-\frac{2\binom{n-2}{k}}{\binom n{k+1}}=1-\frac{2(k+1)(n-k-1)}{n(n-1)}\le1-\frac2n,
\]
where the last inequality follows from $1\le k+1\le n-1$.

On the other hand, for every character that is not large, Lemma~\ref{lem:coefficients} gives
$
\frac{|E_{k+1}(\chi)|}{\binom n{k+1}}\le(n+1)e^{-(k-1)/12},
$
since $\min\{k+1,n-k-1\}\ge k-1$. The Weil bounds give $h\le(\sqrt q+1)^4$ and $n+1\le q+2+4\sqrt q$, so $h(n+1)\le2q^3$ for $q\ge2500$. Since $k\ge60\log q$, the total contribution of these characters to \eqref{eq:genus2-count} is at most
\[
2q^3e^{-(k-1)/12}\le2e^{1/12}q^{-2}<\frac1q\le\frac2n.
\]

Combining all the contribution of the nontrivial characters, the right-hand side of \eqref{eq:genus2-count} is strictly less than $1$. The same conclusion holds if there is no nontrivial large character. This contradiction proves \eqref{eq:n-genus2}. Finally, combining the constructions in Section~\ref{sec:genus2} below, the proof is complete. 
	\end{proof}


	\section{Sharpness of the bounds}\label{sec:construction}
\subsection{Square-root-length NMDS codes on hyperelliptic curves}\label{Sect: sqrt_root_length}
Let $g\ge2$, let $r\ge(2g+2)^2$ be an odd prime power, and put $q=r^2$. Let $X/\F_r$ be the smooth projective curve induced by the affine model $y^2=f(x)$ with $f\in\F_r[x]$ square-free of degree $2g+1$. It is not hard to check that there is a unique point $P_0$ lying above the infinity point of the affine model, which is defined over $\F_r$. 

We still retain the notations in Table~\ref{sec:moments}. Moreover, denote $H=\mathcal J(\F_r)$ as a subgroup of $J$. Recall that $A_j$ denotes the set of effective $\F_q$-rational divisors of degree $j$, as defined before \eqref{Eqn: Wj}. For $0\le j\le g$, let
$
\Theta_j=\{[R-jP_0]:R\in A_j\}\subseteq J.
$
We have
$
\Theta_0\subseteq\Theta_1\subseteq\cdots\subseteq\Theta_g=J,
$
where the last equality follows from Riemann--Roch. Moreover, every degree-zero divisor supported on $X(\F_r)$ has its divisor class in $H$. The following intermediate-coset lemma is the key to our NMDS construction.

\begin{lemma}\label{lem:hyperelliptic-intermediate-coset}
Under the above assumptions, there exists a coset $\mathcal T$ of $H$ in $J$ such that
\[
\mathcal T\cap\Theta_{g-2}=\varnothing,\qquad
\mathcal T\cap\Theta_{g-1}\ne\varnothing.
\]
\end{lemma}

\begin{proof}
For every $b\in J$, we have $b\in H+\Theta_i$ if and only if $(b+H)\cap\Theta_i\ne\varnothing$. It therefore suffices to show that $H+\Theta_{g-2}\subsetneq H+\Theta_{g-1}$.

Suppose that $H+\Theta_{g-2}=H+\Theta_{g-1}$. Since $\Theta_{g-2}+\iota(P)\subseteq\Theta_{g-1}$ for every $P\in X(\F_q)$, our assumption implies that 
$
(H+\Theta_{g-2})+\iota(P)\subseteq H+\Theta_{g-2}.
$
Since translation preserves the size of the sets, hence this inclusion is an equality. Since the point classes $\iota(P)$ generate $J$, it follows that the nonempty set $H+\Theta_{g-2}$ is invariant under translation by every element of $J$. Thus $H+\Theta_{g-2}=J$, and hence $H+\Theta_{g-1}= J$.

On the other hand, by \cite[Lemma~15]{HuMaConstruction}, there exist distinct elements $a_1,a_2\in\F_r$ such that $f(a_1)$ and $f(a_2)$ are nonsquares in $\F_r$. Choose $y_i\in\F_q$ with $y_i^2=f(a_i)$ for $i=1,2$. Then the points $(a_1,y_1)$ and $(a_2,y_2)$ are fixed by the composition of the $r$-power Frobenius and the hyperelliptic involution. Thus \cite[Proposition~14]{HuMaConstruction} yields an element $b\in J$ such that $(b+H)\cap\Theta_{g-1}=\varnothing$. Equivalently,
$
H+\Theta_{g-1}\ne J.
$ This results a contradistinction, hence $H+\Theta_{g-2}\subsetneq H+\Theta_{g-1}$ is true. Thus the proof is done. 
\end{proof}

Let $D=\sum_{P\in X(\F_r)}P$ and $n=\deg D$. By the Hasse--Weil bound and the assumption $\sqrt r\ge2g+2$, we have
$
n<q.
$
Fix an integer $g+1\le k\le n-g$. Choose a coset $\mathcal T$ as in Lemma~\ref{lem:hyperelliptic-intermediate-coset} and an element $w\in\mathcal T\cap\Theta_{g-1}$. Let $E\le D$ be an effective subdivisor of degree $k$. By \cite[Lemma~6]{HuMaConstruction}, since $n<q$ and $k\ge g+1$, there exists an effective $\F_q$-rational divisor $G$ of degree $k+g-1$, with support disjoint from that of $D$, such that
\begin{equation}\label{eq:hyperelliptic-G}
[G-(k+g-1)P_0]=w+[E-kP_0].
\end{equation}
Since $2g-2<\deg G<n$, Riemann--Roch gives $\ell(G)=k$, and the evaluation map is injective. Thus $\CL(D,G)$ has length $n$ and dimension $k$. We now show that it is an NMDS code.

\begin{proposition}\label{prop:hyperelliptic-constructions}
With the above notation, the code $\CL(D,G)$ is an $[n,k,n-k]_q$ NMDS code. In particular, for every fixed genus $g\ge2$ and infinitely many square prime powers $q$, there exist NMDS algebraic geometry codes of length $n=\sqrt q+2gq^{1/4}+1$ for every dimension $g+1\le k\le n-g$.
\end{proposition}

\begin{proof}[Proof of Proposition~\ref{prop:hyperelliptic-constructions}]
By construction, $\CL(D,G)$ has length $n$ and dimension $k$. Since $[G-E-(g-1)P_0]=w\in\Theta_{g-1}$, we have $\ell(G-E)\ge1$. Thus, $d(\CL(D,G))\le n-k$.

On the other hand, let $B$ be any sum of $k+1$ distinct points of $\supp(D)$. Then
$
[G-B-(g-2)P_0]=w+[E-B+P_0]\in\mathcal T.
$
If $\ell(G-B)>0$, then $[G-B-(g-2)P_0]$ belongs to $\Theta_{g-2}$, contradicting the choice of $\mathcal T$. Consequently, $\ell(G-B)=0$ for every such $B$. Consequently, 
$
d(\CL(D,G))=n-k.
$

Since the dimension of $\CL(D,G)^{\perp}$ is $n-k$, to show that this code is NMDS, it suffices to prove that $\CL(D,G)^{\perp}$ has distance $k$. Equivalently, every $k-1$ columns of a generator matrix of $\CL(D,G)$ are linearly independent. We prove this by checking that the evaluation map restricted to any $k-1$ coordinates has rank $k-1$.

Let $B$ be any sum of $k-1$ distinct points of $\supp(D)$. Since $\deg(G-B)=g$, Riemann--Roch gives $\ell(G-B)\ge1$. Notice that our desired conclusion is equivalent to $\ell(G-B)=1$. Suppose that $\ell(G-B)\ge2$. Choose a nonzero function $h\in L(G-B)$ and put $A=(h)+G-B$. Then $A$ is an effective $\F_q$-rational divisor of degree $g$, and $\ell(A)=\ell(G-B)\ge2$. By \cite[Corollary~9]{HuMaConstruction}, $A$ is linearly equivalent to $2P_0+A'$ for some effective $\F_q$-rational divisor $A'$. So we have 
\[
[G-B-gP_0]=[A'-(g-2)P_0]\in\Theta_{g-2}.
\]
However, the defining relation for $G$ also gives
\[
[G-B-gP_0]=w+[E-B-P_0]\in\mathcal T.
\]
This contradicts $\mathcal T\cap\Theta_{g-2}=\varnothing$. Therefore
$
\ell(G-B)=1
$
for every sum $B$ of $k-1$ distinct evaluation points. Thus $d(\CL(D,G)^{\perp})=k$ as claimed. This proves that $\CL(D,G)$ is an NMDS code.

Finally, fix a prime $p\equiv-1\pmod{2(2g+1)}$ and a sufficiently large positive odd integer $m$. Put $r=p^{2m}$ and $q=r^2=p^{4m}$, and consider the curve
\[
X:\quad y^2=1-x^{2g+1}
\]
over $\F_r$. The congruence condition ensures that the defining polynomial is square-free. By \cite[Lemma~19]{HuMaConstruction}, this curve is maximal over $\F_{p^2}$, and, since $m$ is odd, its number of $\F_r$-rational points is
\[
n=\#X(\F_r)=r+1+2gp^m=\sqrt q+2gq^{1/4}+1.
\]
For all sufficiently large odd $m$, we have $r\ge(2g+2)^2$, so the preceding NMDS constructions apply in every dimension $g+1\le k\le n-g$. This completes the proof.
\end{proof}

\subsection{The constructions for genus-two NMDS codes}\label{sec:genus2}
    In this subsection, we construct a family of genus two NMDS codes which achieve the upper bound given in Proposition~\ref{prop:nmds-genus-two}. Thus it completes the proof of this proposition by combining with the proof in Section~\ref{Sect: pf_prop:nmds-genus-two}. 

Let $q=r^2$, where $r$ is an odd prime power satisfying $r\equiv-1\pmod8$. Consider the smooth projective curves $X$ and $Y$ over $\F_q$ induced by the affine equations
\[
X: y^2=u(1-u^4),\qquad Y: v^2=1-t^8.
\]
These curves have genera two and three, respectively. The morphism
\[
\pi:Y\longrightarrow X,\qquad (t,v)\longmapsto(t^2,tv)
\]
is the quotient by the involution $(t,v)\mapsto(-t,-v)$. This involution has no affine fixed points and exchanges the two points at infinity of $Y$, which are rational over $\F_q$. Thus $\pi$ is a geometrically connected unramified double cover, and the unique point $P_0$ at infinity of $X$ splits in $Y$.

Since $8\mid r+1$, taking $t=U^{(r+1)/8}$ and $v=V^{(r+1)/2}$ induces a nonconstant morphism from the Hermitian curve $U^{r+1}+V^{r+1}=1$ to $Y$. The Hermitian curve is maximal over $\F_q$ in the sense that its number of $\F_q$-rational points attains the Hasse--Weil upper bound \cite[p.~138]{GarciaStichtenothXing}. Since maximality descends under nonconstant morphisms of smooth projective curves (see \cite{Stichtenoth} and \cite[p.~138]{GarciaStichtenothXing}), both $X$ and $Y$ are also maximal. Thus we obtain
\[
N:=\#X(\F_q)=q+1+4r,\qquad \#Y(\F_q)=q+1+6r.
\]

Let $A\subseteq X(\F_q)$ be the set of rational points that split completely in via $\pi$, namely, those whose fibers consist of two $\F_q$-rational points. Since $\pi$ is unramified, every point outside $A$ has no rational preimage, we obtain
\[
n:=|A|=\frac{\#Y(\F_q)}2=\frac{\#X(\F_q)}2+\sqrt q.
\]
By Galois theory, or the discussion around \eqref{eq:character-L-series}, there is a nontrivial quadratic character $\chi\in\wh J$ such that $\chi(\iota(P))=1$ if and only if $P\in A$. Since $N-n>0$, we can find $R\in X(\F_q)\setminus A$, for which $\chi(\iota(R))=-1$.

Fix an integer $2\le k\le n-2$, and let $B$ be the sum of $k$ distinct points of $A$. Put $D=\sum_{P\in A}P$.  By weak approximation, there exists an $\F_q$-rational divisor $G\sim B+R$ whose support is disjoint from that of $D$. Since $\deg G=k+1<n$, Riemann--Roch gives $\ell(G)=k$, and the evaluation map is injective. Thus $C=\CL(D,G)$ has dimension $k$.

Since $\widetilde\chi(G)=\widetilde\chi(B+R)=-1$, whereas every sum of $k+1$ distinct evaluation points has character value $1$, no nonzero function in $L(G)$ can vanish at $k+1$ evaluation points or more. Hence $d(C)\ge n-k$. On the other hand, the linear equivalence $G\sim B+R$ gives a nonzero function in $L(G)$ vanishing at exactly the $k$ points of $B$. Consequently, $d(C)=n-k$. By Lemma~\ref{Lem: dual_AG}, the dual code is $\CL(D,G^\perp)$ for a divisor $G^\perp\sim D-G+K_X$ with support disjoint from $D$. Since $X$ is hyperelliptic of genus two, we have $K_X\sim2P_0$. Moreover, $\widetilde\chi(D)=1$ and $\widetilde\chi(G)=-1$, so $\widetilde\chi(G^\perp)=-1$. By the same arguments as above, this excludes any sum of $n-k+1$ distinct evaluation points from being linearly equivalent to $G^\perp$, giving $d(C^\perp)\ge k$. Since $C$ is not MDS, neither is its dual, and therefore $d(C^\perp)=k$. Thus $C$ is NMDS. This completes the proof of the equality assertion in Proposition~\ref{prop:nmds-genus-two}.

\bibliographystyle{plain}
\bibliography{main}

@article {Ball,
    AUTHOR = {Ball, Simeon},
     TITLE = {On sets of vectors of a finite vector space in which every
              subset of basis size is a basis},
   JOURNAL = {J. Eur. Math. Soc. (JEMS)},
  FJOURNAL = {Journal of the European Mathematical Society (JEMS)},
    VOLUME = {14},
      YEAR = {2012},
    NUMBER = {3},
     PAGES = {733--748},
      ISSN = {1435-9855,1435-9863},
   MRCLASS = {15A03 (05B35 51E21 94B05)},
  MRNUMBER = {2911882},
MRREVIEWER = {Gengsheng\ Zhang},
       DOI = {10.4171/JEMS/316},
       URL = {https://doi.org/10.4171/JEMS/316},
}

@article {BDB,
    AUTHOR = {Ball, Simeon and De Beule, Jan},
     TITLE = {On sets of vectors of a finite vector space in which every
              subset of basis size is a basis {II}},
   JOURNAL = {Des. Codes Cryptogr.},
  FJOURNAL = {Designs, Codes and Cryptography. An International Journal},
    VOLUME = {65},
      YEAR = {2012},
    NUMBER = {1-2},
     PAGES = {5--14},
      ISSN = {0925-1022,1573-7586},
   MRCLASS = {15A03 (05B35 51E21 94B05)},
  MRNUMBER = {2943640},
MRREVIEWER = {Gengsheng\ Zhang},
       DOI = {10.1007/s10623-012-9658-6},
       URL = {https://doi.org/10.1007/s10623-012-9658-6},
}

@article {BL,
    AUTHOR = {Ball, Simeon and Lavrauw, Michel},
     TITLE = {Arcs in finite projective spaces},
   JOURNAL = {EMS Surv. Math. Sci.},
  FJOURNAL = {EMS Surveys in Mathematical Sciences},
    VOLUME = {6},
      YEAR = {2019},
    NUMBER = {1-2},
     PAGES = {133--172},
      ISSN = {2308-2151,2308-216X},
   MRCLASS = {51E21 (05B25 94B05)},
  MRNUMBER = {4073887},
MRREVIEWER = {Yue\ Zhou},
       DOI = {10.4171/emss/33},
       URL = {https://doi.org/10.4171/emss/33},
}

@article {Bombieri,
    AUTHOR = {Bombieri, Enrico},
     TITLE = {On exponential sums in finite fields},
   JOURNAL = {Amer. J. Math.},
  FJOURNAL = {American Journal of Mathematics},
    VOLUME = {88},
      YEAR = {1966},
     PAGES = {71--105},
      ISSN = {0002-9327,1080-6377},
   MRCLASS = {14.48 (10.41)},
  MRNUMBER = {200267},
MRREVIEWER = {D.\ J.\ Lewis},
       DOI = {10.2307/2373048},
       URL = {https://doi.org/10.2307/2373048},
}

@article {BTB,
    AUTHOR = {Bruen, A. A. and Thas, J. A. and Blokhuis, A.},
     TITLE = {On {M}.{D}.{S}. codes, arcs in {${\rm PG}(n,q)$} with {$q$}
              even, and a solution of three fundamental problems of {B}.
              {S}egre},
   JOURNAL = {Invent. Math.},
  FJOURNAL = {Inventiones Mathematicae},
    VOLUME = {92},
      YEAR = {1988},
    NUMBER = {3},
     PAGES = {441--459},
      ISSN = {0020-9910,1432-1297},
   MRCLASS = {51E20 (05B25 14G15 62K15)},
  MRNUMBER = {939470},
MRREVIEWER = {David\ G.\ Glynn},
       DOI = {10.1007/BF01393742},
       URL = {https://doi.org/10.1007/BF01393742},
}

@article {HaoChen,
    AUTHOR = {Chen, Hao},
     TITLE = {On the main conjecture of geometric {MDS} codes},
   JOURNAL = {Internat. Math. Res. Notices},
  FJOURNAL = {International Mathematics Research Notices},
    VOLUME = {1994},
      YEAR = {1994},
    NUMBER = {8},
     PAGES = {313--318},
      ISSN = {1073-7928,1687-0247},
   MRCLASS = {94B27 (11T71 14G15)},
  MRNUMBER = {1289576},
MRREVIEWER = {Michael\ M.\ Dediu},
       DOI = {10.1155/S107379289400036X},
       URL = {https://doi.org/10.1155/S107379289400036X},
}

@article {CY,
    AUTHOR = {Chen, Hao and Yau, Stephen S.-T.},
     TITLE = {Contribution to {M}unuera's problem on the main conjecture of
              geometric hyperelliptic {MDS} codes},
   JOURNAL = {IEEE Trans. Inform. Theory},
  FJOURNAL = {Institute of Electrical and Electronics Engineers.
              Transactions on Information Theory},
    VOLUME = {43},
      YEAR = {1997},
    NUMBER = {4},
     PAGES = {1349--1354},
      ISSN = {0018-9448,1557-9654},
   MRCLASS = {94B27},
  MRNUMBER = {1454971},
       DOI = {10.1109/18.605607},
       URL = {https://doi.org/10.1109/18.605607},
}

@incollection {deBoer,
    AUTHOR = {de Boer, Mario A.},
     TITLE = {M{DS} codes from hyperelliptic curves},
 BOOKTITLE = {Arithmetic, geometry and coding theory ({L}uminy, 1993)},
     PAGES = {23--34},
 PUBLISHER = {de Gruyter, Berlin},
      YEAR = {1996},
      ISBN = {3-11-014616-9},
   MRCLASS = {11G20 (11T71 14G15 94B27)},
  MRNUMBER = {1394923},
MRREVIEWER = {Marc\ Perret},
}

@article {DL,
    AUTHOR = {Dodunekov, Stefan and Landgev, Ivan},
     TITLE = {On near {MDS} codes},
   JOURNAL = {J. Geom.},
  FJOURNAL = {Journal of Geometry},
    VOLUME = {54},
      YEAR = {1995},
    NUMBER = {1-2},
     PAGES = {30--43},
      ISSN = {0047-2468,1420-8997},
   MRCLASS = {94B05 (05B25 51E22)},
  MRNUMBER = {1358272},
MRREVIEWER = {E.\ F.\ Assmus, Jr.},
       DOI = {10.1007/BF01222850},
       URL = {https://doi.org/10.1007/BF01222850},
}

@article {DLsmall,
    AUTHOR = {Dodunekov, S. M. and Landjev, I. N.},
     TITLE = {Near-{MDS} codes over some small fields},
   JOURNAL = {Discrete Math.},
  FJOURNAL = {Discrete Mathematics},
    VOLUME = {213},
      YEAR = {2000},
    NUMBER = {1-3},
     PAGES = {55--65},
      ISSN = {0012-365X,1872-681X},
   MRCLASS = {94B60},
  MRNUMBER = {1755411},
       DOI = {10.1016/S0012-365X(99)00168-5},
       URL = {https://doi.org/10.1016/S0012-365X(99)00168-5},
      NOTE = {Selected topics in discrete mathematics (Warsaw, 1996)},
}

@article {Goppa,
    AUTHOR = {Goppa, V. D.},
     TITLE = {Codes on algebraic curves},
   JOURNAL = {Dokl. Akad. Nauk SSSR},
  FJOURNAL = {Doklady Akademii Nauk SSSR},
    VOLUME = {259},
      YEAR = {1981},
    NUMBER = {6},
     PAGES = {1289--1290},
      ISSN = {0002-3264},
   MRCLASS = {94B05 (14H45)},
  MRNUMBER = {628795},
MRREVIEWER = {E.\ J. F. Primrose},
}

@article {GS,
    AUTHOR = {Garc\'{\i}a, Arnaldo and Stichtenoth, Henning},
     TITLE = {A tower of {A}rtin-{S}chreier extensions of function fields
              attaining the {D}rinfel'd-{V}l\u{a}du\c{t} bound},
   JOURNAL = {Invent. Math.},
  FJOURNAL = {Inventiones Mathematicae},
    VOLUME = {121},
      YEAR = {1995},
    NUMBER = {1},
     PAGES = {211--222},
      ISSN = {0020-9910,1432-1297},
   MRCLASS = {11G20 (11R58 14G15 94B27)},
  MRNUMBER = {1345289},
MRREVIEWER = {Jos\'{e}\ Felipe\ Voloch},
       DOI = {10.1007/BF01884295},
       URL = {https://doi.org/10.1007/BF01884295},
}

@article {HR,
    AUTHOR = {Han, Dongchun and Ren, Yuan},
     TITLE = {The maximal length of {$q$}-ary {MDS} elliptic codes is close
              to {$\frac q2$}},
   JOURNAL = {Int. Math. Res. Not. IMRN},
  FJOURNAL = {International Mathematics Research Notices. IMRN},
    VOLUME = {2024},
      YEAR = {2024},
    NUMBER = {11},
     PAGES = {9036--9043},
      ISSN = {1073-7928,1687-0247},
   MRCLASS = {14G50 (05E14 11B30 11T71 94B05)},
  MRNUMBER = {4756098},
MRREVIEWER = {Juan\ Tena},
       DOI = {10.1093/imrn/rnad271},
       URL = {https://doi.org/10.1093/imrn/rnad271},
}

@article {HRTight,
    AUTHOR = {Han, Dongchun and Ren, Yuan},
     TITLE = {A tight upper bound for the maximal length of {MDS} elliptic
              codes},
   JOURNAL = {IEEE Trans. Inform. Theory},
  FJOURNAL = {Institute of Electrical and Electronics Engineers.
              Transactions on Information Theory},
    VOLUME = {69},
      YEAR = {2023},
    NUMBER = {2},
     PAGES = {819--822},
      ISSN = {0018-9448,1557-9654},
   MRCLASS = {94B05},
  MRNUMBER = {4564628},
MRREVIEWER = {Yan\ Liu},
       DOI = {10.1109/tit.2022.3211543},
       URL = {https://doi.org/10.1109/tit.2022.3211543},
}

@article {HK,
    AUTHOR = {Hirschfeld, J. W. P. and Korchm\'{a}ros, G.},
     TITLE = {On the number of rational points on an algebraic curve over a
              finite field},
   JOURNAL = {Bull. Belg. Math. Soc. Simon Stevin},
  FJOURNAL = {Bulletin of the Belgian Mathematical Society. Simon Stevin},
    VOLUME = {5},
      YEAR = {1998},
    NUMBER = {2-3},
     PAGES = {313--340},
      ISSN = {1370-1444,2034-1970},
   MRCLASS = {11G20 (14G05 14G15 51E15 51E20)},
  MRNUMBER = {1630034},
MRREVIEWER = {Chao\ Ping\ Xing},
       URL = {http://projecteuclid.org/euclid.bbms/1103409013},
      NOTE = {Finite geometry and combinatorics (Deinze, 1997)},
}

@book {Liu,
    AUTHOR = {Liu, Qing},
     TITLE = {Algebraic geometry and arithmetic curves},
    SERIES = {Oxford Graduate Texts in Mathematics},
    VOLUME = {6},
      NOTE = {Translated from the French by Reinie Ern\'{e},
              Oxford Science Publications},
 PUBLISHER = {Oxford University Press, Oxford},
      YEAR = {2002},
     PAGES = {xvi+576},
      ISBN = {0-19-850284-2},
   MRCLASS = {14-01 (11G30 14A05 14A15 14Gxx 14Hxx)},
  MRNUMBER = {1917232},
MRREVIEWER = {C\'{\i}cero\ Carvalho},
}

@book {MS,
    AUTHOR = {MacWilliams, F. J. and Sloane, N. J. A.},
     TITLE = {The theory of error-correcting codes. {I}},
    SERIES = {North-Holland Mathematical Library, Vol. 16},
 PUBLISHER = {North-Holland Publishing Co., Amsterdam-New York-Oxford},
      YEAR = {1977},
     PAGES = {i--xv and 1--369},
      ISBN = {0-444-85009-0},
   MRCLASS = {94A10},
  MRNUMBER = {465509},
MRREVIEWER = {Ian\ Blake},
}

@article {Munuera,
    AUTHOR = {Munuera, Carlos},
     TITLE = {On the main conjecture on geometric {MDS} codes},
   JOURNAL = {IEEE Trans. Inform. Theory},
  FJOURNAL = {Institute of Electrical and Electronics Engineers.
              Transactions on Information Theory},
    VOLUME = {38},
      YEAR = {1992},
    NUMBER = {5},
     PAGES = {1573--1577},
      ISSN = {0018-9448,1557-9654},
   MRCLASS = {94B05 (11T71 14H52 94B27)},
  MRNUMBER = {1178199},
MRREVIEWER = {Francesco\ Fabris},
       DOI = {10.1109/18.149510},
       URL = {https://doi.org/10.1109/18.149510},
}

@article {MunueraEll,
    AUTHOR = {Munuera, Carlos},
     TITLE = {On {MDS} elliptic codes},
   JOURNAL = {Discrete Math.},
  FJOURNAL = {Discrete Mathematics},
    VOLUME = {117},
      YEAR = {1993},
    NUMBER = {1-3},
     PAGES = {279--286},
      ISSN = {0012-365X,1872-681X},
   MRCLASS = {94B27 (11T71)},
  MRNUMBER = {1226151},
       DOI = {10.1016/0012-365X(93)90344-S},
       URL = {https://doi.org/10.1016/0012-365X(93)90344-S},
}

@article {Segre,
    AUTHOR = {Segre, Beniamino},
     TITLE = {Ovals in a finite projective plane},
   JOURNAL = {Canadian J. Math.},
  FJOURNAL = {Canadian Journal of Mathematics. Journal Canadien de
              Math\'{e}matiques},
    VOLUME = {7},
      YEAR = {1955},
     PAGES = {414--416},
      ISSN = {0008-414X,1496-4279},
   MRCLASS = {48.0X},
  MRNUMBER = {71034},
MRREVIEWER = {Marshall\ Hall, Jr.},
       DOI = {10.4153/CJM-1955-045-x},
       URL = {https://doi.org/10.4153/CJM-1955-045-x},
}

@book {SerreClassFields,
    AUTHOR = {Serre, Jean-Pierre},
     TITLE = {Algebraic groups and class fields},
    SERIES = {Graduate Texts in Mathematics},
    VOLUME = {117},
      NOTE = {Translated from the French},
 PUBLISHER = {Springer-Verlag, New York},
      YEAR = {1988},
     PAGES = {x+207},
      ISBN = {0-387-96648-X},
   MRCLASS = {14K30 (11G45 11S31 14L40)},
  MRNUMBER = {918564},
       DOI = {10.1007/978-1-4612-1035-1},
       URL = {https://doi.org/10.1007/978-1-4612-1035-1},
}

@article {Shokrollahi,
    AUTHOR = {Shokrollahi, Mohammad Amin},
     TITLE = {Minimum distance of elliptic codes},
   JOURNAL = {Adv. Math.},
  FJOURNAL = {Advances in Mathematics},
    VOLUME = {93},
      YEAR = {1992},
    NUMBER = {2},
     PAGES = {251--281},
      ISSN = {0001-8708,1090-2082},
   MRCLASS = {11T71 (11G25)},
  MRNUMBER = {1164709},
MRREVIEWER = {\'{A}.\ G.\ Horv\'{a}th},
       DOI = {10.1016/0001-8708(92)90029-K},
       URL = {https://doi.org/10.1016/0001-8708(92)90029-K},
}

@book {Stichtenoth,
    AUTHOR = {Stichtenoth, Henning},
     TITLE = {Algebraic function fields and codes},
    SERIES = {Graduate Texts in Mathematics},
    VOLUME = {254},
   EDITION = {Second},
 PUBLISHER = {Springer-Verlag, Berlin},
      YEAR = {2009},
     PAGES = {xiv+355},
      ISBN = {978-3-540-76877-7},
   MRCLASS = {14H05 (11R58 11T71 14G15 14G50 94B27)},
  MRNUMBER = {2464941},
}

@article {ST,
    AUTHOR = {Storme, L. and Thas, J. A.},
     TITLE = {M.{D}.{S}. codes and arcs in {${\rm PG}(n,q)$} with {$q$}
              even: an improvement of the bounds of {B}ruen, {T}has, and
              {B}lokhuis},
   JOURNAL = {J. Combin. Theory Ser. A},
  FJOURNAL = {Journal of Combinatorial Theory. Series A},
    VOLUME = {62},
      YEAR = {1993},
    NUMBER = {1},
     PAGES = {139--154},
      ISSN = {0097-3165,1096-0899},
   MRCLASS = {51E21 (51E22 94B27 94B65)},
  MRNUMBER = {1198385},
MRREVIEWER = {J.\ W. P. Hirschfeld},
       DOI = {10.1016/0097-3165(93)90076-K},
       URL = {https://doi.org/10.1016/0097-3165(93)90076-K},
}

@article {Thas1975,
    AUTHOR = {Thas, J. A.},
     TITLE = {Some results concerning {$\{(q+1)(n-1);\,n\}$}-arcs and
              {$\{(q+1)$} {$(n-1)+1;\,n\}$}-arcs in finite projective planes
              of order {$q$}},
   JOURNAL = {J. Combinatorial Theory Ser. A},
  FJOURNAL = {Journal of Combinatorial Theory. Series A},
    VOLUME = {19},
      YEAR = {1975},
    NUMBER = {2},
     PAGES = {228--232},
      ISSN = {0097-3165},
   MRCLASS = {50D30},
  MRNUMBER = {377682},
MRREVIEWER = {W.\ M.\ Kantor},
       DOI = {10.1016/s0097-3165(75)80012-4},
       URL = {https://doi.org/10.1016/s0097-3165(75)80012-4},
}

@book {TV,
    AUTHOR = {Tsfasman, M. A. and Vl\u{a}du\c{t}, S. G.},
     TITLE = {Algebraic-geometric codes},
    SERIES = {Mathematics and its Applications (Soviet Series)},
    VOLUME = {58},
      NOTE = {Translated from the Russian by the authors},
 PUBLISHER = {Kluwer Academic Publishers Group, Dordrecht},
      YEAR = {1991},
     PAGES = {xxiv+667},
      ISBN = {0-7923-0727-5},
   MRCLASS = {94B27 (11T71 94-02)},
  MRNUMBER = {1186841},
MRREVIEWER = {Simon\ N.\ Litsyn},
       DOI = {10.1007/978-94-011-3810-9},
       URL = {https://doi.org/10.1007/978-94-011-3810-9},
}

@article {TVZ,
    AUTHOR = {Tsfasman, M. A. and Vl\u{a}du\c{t}, S. G. and Zink, Th.},
     TITLE = {Modular curves, {S}himura curves, and {G}oppa codes, better
              than {V}arshamov-{G}ilbert bound},
   JOURNAL = {Math. Nachr.},
  FJOURNAL = {Mathematische Nachrichten},
    VOLUME = {109},
      YEAR = {1982},
     PAGES = {21--28},
      ISSN = {0025-584X,1522-2616},
   MRCLASS = {11T71 (14G15 14H25 94B05)},
  MRNUMBER = {705893},
MRREVIEWER = {Yasutaka\ Ihara},
       DOI = {10.1002/mana.19821090103},
       URL = {https://doi.org/10.1002/mana.19821090103},
}

@article {Walker,
    AUTHOR = {Walker, Judy L.},
     TITLE = {A new approach to the main conjecture on algebraic-geometric
              {MDS} codes},
   JOURNAL = {Des. Codes Cryptogr.},
  FJOURNAL = {Designs, Codes and Cryptography. An International Journal},
    VOLUME = {9},
      YEAR = {1996},
    NUMBER = {1},
     PAGES = {115--120},
      ISSN = {0925-1022,1573-7586},
   MRCLASS = {94B27 (11G20 11T71)},
  MRNUMBER = {1412183},
       DOI = {10.1023/A:1027358511882},
       URL = {https://doi.org/10.1023/A:1027358511882},
      NOTE = {Second Upper Michigan Combinatorics Workshop on Designs, Codes
              and Geometries (Houghton, MI, 1994)},
}

@article {ZhuEtAl,
    AUTHOR = {Zhu, Yunlong and Chen, Hao and Peng, Fan and Zhang, Fangguo
              and Zhao, Chang-An},
     TITLE = {New explicit constructions of {MDS} and {NMDS} algebraic
              geometry codes},
   JOURNAL = {IEEE Trans. Inform. Theory},
  FJOURNAL = {Institute of Electrical and Electronics Engineers.
              Transactions on Information Theory},
    VOLUME = {72},
      YEAR = {2026},
    NUMBER = {9},
     PAGES = {6450--6465},
      ISSN = {0018-9448,1557-9654},
   MRCLASS = {94B27},
  MRNUMBER = {5119224},
       DOI = {10.1109/tit.2026.3712030},
       URL = {https://doi.org/10.1109/tit.2026.3712030},
}

@article {GarciaStichtenothXing,
    AUTHOR = {Garcia, Arnaldo and Stichtenoth, Henning and Xing, Chao-Ping},
     TITLE = {On subfields of the {H}ermitian function field},
   JOURNAL = {Compositio Math.},
  FJOURNAL = {Compositio Mathematica},
    VOLUME = {120},
      YEAR = {2000},
    NUMBER = {2},
     PAGES = {137--170},
      ISSN = {0010-437X,1570-5846},
   MRCLASS = {11R58 (11G20 14G15)},
  MRNUMBER = {1739176},
MRREVIEWER = {Gabriel\ D.\ Villa-Salvador},
       DOI = {10.1023/A:1001736016924},
       URL = {https://doi.org/10.1023/A:1001736016924},
}

@misc{CHMDS,
    AUTHOR = {Chen, Bocong and Huang, Jing},
     TITLE = {On the {MDS} conjecture in even characteristic},
      YEAR = {2026},
HOWPUBLISHED = {arXiv preprint arXiv:2610.28208},
       DOI = {10.48550/arXiv.2610.28208},
       URL = {https://arxiv.org/abs/2610.28208},
}

@misc{CHRestricted,
    AUTHOR = {Chen, Bocong and Huang, Jing},
     TITLE = {Critical numbers for restricted sumsets: rigidity and collapse
              in finite abelian groups},
      YEAR = {2026},
HOWPUBLISHED = {arXiv preprint arXiv:2602.10402},
       DOI = {10.48550/arXiv.2602.10402},
       URL = {https://arxiv.org/abs/2602.10402},
}

@misc{Chen2,
    AUTHOR = {Chen, Haojie},
     TITLE = {A uniform upper bound on the length of {MDS} codes from
              genus-two curves},
      YEAR = {2026},
HOWPUBLISHED = {arXiv preprint arXiv:2610.39364},
       DOI = {10.48550/arXiv.2610.39364},
       URL = {https://arxiv.org/abs/2610.39364},
}

@misc{CHM,
    AUTHOR = {Chen, Haojie and Hu, Chuangqiang and Ma, Guorui},
     TITLE = {Length bounds for {MDS} algebraic geometry codes of arbitrary
              positive genus},
      YEAR = {2026},
HOWPUBLISHED = {arXiv preprint arXiv:2610.39434},
       DOI = {10.48550/arXiv.2610.39434},
       URL = {https://arxiv.org/abs/2610.39434},
}

@misc{Fan,
    AUTHOR = {Fan, Xiang},
     TITLE = {A linear-in-{$q$} range of dimensions for the {MDS} conjecture
              over {$\mathbb F_q$} in odd characteristic},
      YEAR = {2026},
HOWPUBLISHED = {arXiv preprint arXiv:2610.03740},
       DOI = {10.48550/arXiv.2610.03740},
       URL = {https://arxiv.org/abs/2610.03740v2},
}

@misc{HRCover,
    AUTHOR = {Han, Dongchun and Ren, Yuan},
     TITLE = {When do restricted sumsets cover a finite abelian group?},
      YEAR = {2026},
HOWPUBLISHED = {arXiv preprint arXiv:2610.14496},
       DOI = {10.48550/arXiv.2610.14496},
       URL = {https://arxiv.org/abs/2610.14496},
}

@misc{HuMaConstruction,
    AUTHOR = {Hu, Chuangqiang and Ma, Guorui},
     TITLE = {{MDS} algebraic geometry codes of square-root length on
              hyperelliptic curves},
      YEAR = {2026},
HOWPUBLISHED = {arXiv preprint arXiv:2611.01546},
       DOI = {10.48550/arXiv.2611.01546},
       URL = {https://arxiv.org/abs/2611.01546},
}

@inproceedings{LWZ,
    AUTHOR = {Li, Jiyou and Wan, Daqing and Zhang, Jun},
     TITLE = {On the minimum distance of elliptic curve codes},
 BOOKTITLE = {Proc. IEEE Int. Symp. Information Theory (ISIT)},
      YEAR = {2015},
     PAGES = {2391--2395},
 PUBLISHER = {IEEE},
       DOI = {10.1109/ISIT.2015.7282884},
       URL = {https://arxiv.org/abs/1501.01138v2},
}

@article{KT,
    AUTHOR = {Katsman, G. L. and Tsfasman, M. A.},
     TITLE = {Spectra of algebraic-geometric codes},
   JOURNAL = {Problems of Information Transmission},
    VOLUME = {23},
      YEAR = {1987},
     PAGES = {262--275},
       URL = {https://www.mathnet.ru/eng/ppi827},
}

@misc{Kumar,
    AUTHOR = {Kumar, Pranava},
     TITLE = {The {MDS} conjecture over finite fields of even order},
      NOTE = {Palomar Registry, PALOMAR-2026-10-15-000014},
      YEAR = {2026},
       URL = {https://palomar-registry.org/entry?id=PALOMAR-2026-10-15-000014},
}

@article{LR,
    AUTHOR = {Landjev, I. and Rousseva, A.},
     TITLE = {On the maximal length of near-{MDS} codes},
   JOURNAL = {Mathematics and Education in Mathematics},
    VOLUME = {55},
      YEAR = {2026},
     PAGES = {405--413},
       DOI = {10.55630/mem.2026.55.405-413},
       URL = {https://doi.org/10.55630/mem.2026.55.405-413},
}

@misc{XRH,
      AUTHOR = {Xia, Lingjie and Ren, Yuan and Han, Dongchun},
       TITLE = {{MDS} and {NMDS} Main Conjectures for {AG} Codes},
        YEAR = {2026},
HOWPUBLISHED = {Unpublished manuscript},
}

\end{document}